\documentclass[12pt, reqno]{amsart}

\usepackage{amssymb,latexsym,amsmath,amsfonts}
\usepackage{amsmath,amscd,amsthm,amssymb}

\numberwithin{equation}{section}

\newcommand{\Ff}{{\mathbb F}}
\newcommand{\Cc}{{\mathbb C}}

\newcommand{\p}{{\mathbf P}}
\newcommand{\q}{{\mathbf Q}}
\newcommand{\cc}{{\mathcal C}}
\newcommand{\M}{{\mathbf M}}
\newcommand{\A}{{\mathbf A}}

\newcommand{\C}{{\mathbf C}}
\newcommand{\D}{{\mathbf D}}

\newcommand{\ii}{{\mathbf I}}
\newcommand{\jj}{{\mathbf J}}

 \newcommand{\PG}{{\mathrm{PG}}}
 
 \newcommand{\GL}{{\mathrm{GL}}}
 \newcommand{\PSL}{{\mathrm{PSL}}}
 \newcommand{\PGL}{{\mathrm{PGL}}}

\DeclareMathOperator{\Ker}{Ker}
\DeclareMathOperator{\Ima}{Im}

\DeclareMathOperator{\Ind}{Ind}

\newcommand{\Proof}{ \noindent{\bf Proof:}\quad }

\def\QED{\qed\medskip\par}
\newtheorem{Theorem} {Theorem} [section]
\newtheorem{Proposition} [Theorem] {Proposition}
\newtheorem{Lemma} [Theorem] {Lemma}
\newtheorem{Corollary} [Theorem] {Corollary}
\newtheorem{Conjecture}[Theorem]{Conjecture}

\newtheorem{Definition}[Theorem]{Definition}
\newtheorem{Remark} [Theorem] {Remark}

\def\Sq{\;\Box}
\def\Nsq{\;\not\!\!\Box}

\begin{document}

\title[Dimensions of Binary Codes]
{On Binary Codes from Conics in $\PG(2,q)$}

\author[Madison and Wu]{Adonus L. Madison and Junhua Wu$^{{\dagger,\star}}$}

\thanks{$^{\dagger}$Research supported in part by NSF 
HBCU-UP Grant Award $\#0929257$ at Lane College.}

\thanks{$^{\star}$Corresponding author.}

\address{Department of Mathematics, Lane College, Jackson, TN, USA}
\email{jwu@lanecollege.edu}

\address{Lane College, Jackson, TN, USA}
\email{adonus$\_$madison@lanecollege.edu}


\begin{abstract}
Let $\A$ be the incidence matrix of passant lines and internal 
points with respect to a conic in $\PG(2,q)$, where $q$ is an 
odd prime power. In this article, we study both geometric and 
algebraic properties of the column $\Ff_2$-null space 
$\mathcal{L}$ of $\A$. In particular, using methods from both 
finite geometry and modular presentation theory, we manage
to compute the dimension of $\mathcal{L}$, which provides
a proof for the conjecture on the dimension of the binary 
code generated by $\mathcal{L}$. 

\end{abstract}

\maketitle

\section{Introduction}\label{intro}

Let $\PG(2,q)$ be the classical projective plane of order $q$ 
with underlying $3$-dimensional vector space $V$ over $\Ff_q$, 
the finite field of order $q$. Throughout this article, 
$\PG(2,q)$ is represented via homogeneous coordinates. Namely, 
a point is written as a non-zero vector $(a_0,a_1,a_2)$ and 
a line is written as $[b_0,b_1,b_2]$ where not all $b_i$ 
($i=1,2,3$) are zero. The set of points
\begin{equation}\label{conic}
\mathcal{O}:=\{(1,r,r^2)\mid r\in \Ff_q\}\cup\{(0,0,1)\}
\end{equation}
give rise to a geometric object called {\it conic} in 
$\PG(2,q)$ \cite{hir}. The above set also comprises the 
projective solutions of the nondegenarate quadratic equation
\begin{equation}\label{quadric}
Q(X_0, X_1, X_2)=X_1^2-X_0 X_2
\end{equation}
over $\Ff_q$. With respect to $\mathcal{O}$, the lines of 
$\PG(2,q)$ are partitioned into passant lines ($Pa$), tangent 
lines ($T$), and secant lines ($Se$) accordingly as the sizes of 
their intersections with $\mathcal{O}$ are $0$, $1$, or $2$. 
Similarly, points are partitioned into internal points 
($I$), conic points ($\mathcal{O}$), and external points ($E$) 
accordingly as the numbers of tangent lines on which they 
lie are $0$, $1$, or $2$.

In \cite{keith1}, one low-density parity-check binary code was 
constructed using the column $\Ff_2$-null space $\mathcal{L}$ 
of the incidence matrix $\A$ of passant lines and internal 
points with respect to $\mathcal{O}$.
With the help of computer software Magma, the authors made a 
conjecture on the dimension of $\mathcal{L}$ as follows:
\begin{Conjecture}\cite[Conjecture 4.7]{keith1}\label{conj}
Let $\mathcal{L}$ be the $\Ff_2$-null space of $\A$, and let 
$\dim_{\Ff_2}(\mathcal{L})$ be the dimension of $\mathcal{L}$. 
Then
$$\dim_{\Ff_2}(\mathcal{L})=\frac{(q-1)^2}{4}.$$
\end{Conjecture}

The goal of this article is to confirm Cojecture~\ref{conj}. 
Apart from the above conjecture, the dimensions of the column 
$\Ff_2$-null spaces of the incidence matrices of external 
points versus secant lines, external points versus passant 
lines, and passant lines versus external points, were 
conjectured in the aforementioned paper \cite{keith1}, and 
have been established in \cite{swx} and \cite{wu}, respectively. 
To start, we recall that the automorphism group $G$ of 
$\mathcal{O}$ is isomorphic to $\PGL(2,q)$, and that the 
normal subgroup $H$ of $G$ is isomorphic to $\PSL(2,q)$.

Let $F$ be an algebraic closure of $\Ff_2$.
The idea of proving Conjecture~\ref{conj} is to first 
realize $\mathcal{L}$ as an $FH$-module and then decompose 
it into a direct sum of its certain submodules whose 
dimensions can be obtained easily. More concretely 
speaking, 
we view $\A$ as the matrix of the following 
homomorphism $\phi$ of free $F$-modules:
$$\phi: F^I\rightarrow F^I$$
which first sends an internal point to the formal 
sum of all internal points on its polar, and then 
extends linearly to the whole of $F^I$. Additionally, 
it can be shown that $\phi$ is indeed an $FH$-module. 
Consequently, computing the dimension of the column 
$\Ff_2$-null space amounts to finding the dimension of 
the $F$-null space of $\phi$. 
To this end, we investigate the underlying $FH$-module 
structure of $\mathcal{L}$ by applying Brauer's theory
on the $2$-blocks of $H$ and arrive at a convenient 
decomposition of $\mathcal{L}$.

This article is organized in the following way. In Section 2, 
we establish that the matrix $\A$ satisfying the relation 
$\A^3\equiv \A \pmod 2$ under certain orderings of its rows 
and columns; this relation, in turn, reveals a 
geometric discription of $\Ker(\phi)$ as well as yields a set 
of generating elements of $\Ker(\phi)$ in terms of the concept 
of internal neighbors. In Section 3, the parity of intersection 
sizes of certain subsets of $H$ with the conjugacy classes of 
$H$ are computed. We will then in Section 4 review several facts
about the $2$-blocks of $\PSL(2,q)$ and the block idempotents 
of the $2$-blocks of $\PSL(2,q)$; the detailed calculations of 
the $2$-block idempotents were performed in \cite{swx}. Combining 
the results in Sections 3 and 4 with Brauer's theory on blocks, 
we are able to decompose $\Ker(\phi)$ into a direct sum of all 
non-isomorphic simple $FH$-modules or this sum plus a trivial 
module depending on $q$. Consequently, the dimension 
of $\mathcal{L}$ follows as a lemma.

\section{Geometry of Conics}

First we recall several well-known results related to the 
geometry of conics in $\PG(2,q)$ with $q$ odd. The books 
\cite{hp} and \cite{hir} are the references for what follows.

A {\it collineation} of $\PG(2,q)$ is an automorphism of
$\PG(2,q)$, which is a bijection from the set of all points 
and all lines of $\PG(2,q)$  to itself that maps a point to 
a line and a line to a point, and preserves incidence. 
It is well known that each element of $\GL(3,q)$, the group 
of all $3\times 3$ non-singular matrices over $\Ff_q$, 
induces a collineation of $\PG(2,q)$. The proof of the 
following lemma is straightforward.

\begin{Lemma}\label{action1} Let
$\p=(a_0,a_1,a_2)$ and $\ell=[b_0,b_1,b_2]$ be a point and a
line of $\PG(2,q)$, respectively. Suppose that $\theta$ is a 
collineation of $\PG(2,q)$ that is induced by 
$\mathbf{D}\in \GL(3,q)$. If we use $\p^\theta$ and $\ell^\theta$ 
to denote the images of $\p$ and $\ell$ under $\theta$, 
respectively, then
$\p^\theta = (a_0,a_1,a_2)^\theta = (a_0,a_1,a_2)\mathbf{D}$
and
$\ell^\theta = [b_0,b_1,b_2]^\theta =[c_0,c_1,c_2]$,
where $c_0,c_1,c_2$ correspond to the first, the second, and 
the third coordinate of the vector 
$\mathbf{D}^{-1}(b_0,b_1,b_2)^\top$, respectively.
\end{Lemma}

A {\it correlation} of $\PG(2,q)$ is a bijection from the set of
points to the set of lines as well as the set of lines to the set
of points that reverses inclusion. A {\it polarity}
of $\PG(2,q)$ is a correlation of order $2$. The image of a point
$\p$ under a correlation $\sigma$ is denoted by $\p^\sigma$, and
that of a line $\ell$ is denoted by $\ell^\sigma$. It can be shown
\cite[p.~181]{hir} that the non-degenerate quadratic form
$Q(X_0,X_1,X_2)$ = $X_1^2-X_0X_2$ induces a polarity $\sigma$ 
(or $\perp$) of $\PG(2,q)$, which can be represented by the matrix
\begin{equation}\label{matrix_M}
\M=\left(\begin{array}{rrr}
0 & 0 & -\frac{1}{2} \\
0 & 1 & 0 \\
-\frac{1}{2} & 0 & 0 \\
\end{array}\right).
\end{equation}

\begin{Lemma}\label{lemB}$(${\rm{\cite[p.~47]{hp}}}$)$
Let $\p=(a_0,a_1,a_2)$ and $\ell=[b_0,b_1,b_2])$ be a point and a
line of $\PG(2,q)$, respetively. If $\sigma$ is the polarity 
represented by the non-singular symmetric matrix $\M$ in 
$($\ref{matrix_M}$)$, then
$\p^{\sigma} = (a_0,a_1,a_2)^{\sigma} =[c_0,c_1,c_2]$
and
$\ell^{\sigma} = [b_0,b_1,b_2]^{\sigma} =(b_0,b_1,b_2)\M^{-1}$,
where $c_0,c_1,c_2$ correspond to the first, the second, the third 
coordinate of the column vector $\M(a_0,a_1,a_2)^\top$, respectively.
\end{Lemma}
For example, if $\p=(x,y,z)$ is a point of $\PG(2,q)$,
then its image under $\sigma$ is $\p^{\sigma}=[z,-2y,x]$.

For convenience, we will denote the set of all non-zero squares of
$\Ff_q$ by $\Sq_q$, and the set of non-squares by $\Nsq_q$; also,
$\Ff_q^*$ is the set of non-zero elements of $\Ff_q$.

\begin{Lemma}\label{bijection}$(${\rm{\cite[p.~181--182]{hir}}}$)$
Assume that $q$ is odd.
\begin{itemize}
\item[(i)] The polarity $\sigma$ above defines three bijections; 
that is,
$\sigma:\;I\rightarrow\;Pa$, $\sigma:\;E\rightarrow\;Se$, and
$\sigma:\;\mathcal{O}\rightarrow\;T$ are all bijections.
\item[(ii)] A line $[b_0,b_1,b_2]$ of $\PG(2,q)$ is a passant,
a tangent, or a secant to $\mathcal{O}$ if and only if 
$b_1^2-4b_0b_2 \in\Nsq_q$, $b_1^2-4b_0b_2 = 0$, or 
$b_1^2-4b_0b_2\in \Sq_q$, respectively.
\item[(iii)] A point $(a_0,a_1,a_2)$ of $\PG(2,q)$ is internal,
absolute, or external if and only if $a_1^2-a_0a_2 \in \Nsq_q$, 
$a_1^2-a_0a_2 =0$, or $a_1^2-a_0a_2 \in \Sq_q$, respectively.
\end{itemize}
\end{Lemma}

The results in the following lemma can be obtained by simple
counting; see {\rm \cite{hir}} for more details and related 
results.
\begin{Lemma}{\rm (\cite[p.~170]{hir})}
Using the above notation, we have
\begin{equation}\label{number}
|T|=|\mathcal{O}|=q+1,\;|Pa|=|I|=\frac{q(q-1)}{2},\;\text{and}\;
|Se|=|E|=\frac{q(q+1)}{2}.
\end{equation}
Also, we have the following tables:
\begin{table}[htp]
\begin{center}
\caption{Number of points on lines of various types}
\bigskip
\begin{tabular}{cccc}
\hline
{Name} & {Absolute points} & {External points} & {Internal points} \\
\hline
{Tangent lines} & $1$ & $q$ & $0$ \\
{Secant lines} & $2$ & $\frac{q-1}{2}$ & $\frac{q-1}{2}$ \\
{Passant lines} & $0$ & $\frac{q+1}{2}$ & $\frac{q+1}{2}$\\
\hline
\end{tabular}
\label{tab1}
\end{center}
\end{table}

\begin{table}[htp]
\begin{center}
\caption{Number of lines through points of various types}
\bigskip
\begin{tabular}{cccc}
\hline
{Name} & {Tangent lines} & {Secant lines} & {Skew lines} \\
\hline
{Absolute points} & $1$ & $q$ & $0$ \\
{External points} & $2$ & $\frac{q-1}{2}$ & $\frac{q-1}{2}$ \\
{Internal points} & $0$ & $\frac{q+1}{2}$ & $\frac{q+1}{2}$\\
\hline
\end{tabular}
\label{tab2}
\end{center}
\end{table}
\end{Lemma}

\subsection{The incidence matrix}

Let $G$ be the automorphism group of $\mathcal{O}$ in $\PGL(3,q)$
(i.e. the subgroup of $\PGL(3,q)$ fixing ${\mathcal O}$ setwise). 

\begin{Lemma}\cite[p. 158]{hir}
$G\cong \PGL(2,q)$.
\end{Lemma}

We define
\begin{equation}\label{grouph}
H:=
\left.\left\{\left(\begin{array}{ccc}
a^2 & ab & b^2\\
2ac & ad+bc & 2bd \\
c^2& cd & d^2\\
 \end{array}\right)\right| a, b, c, d \in \Ff_q, ad-bc=1\right\}.
\end{equation}
In the rest of the article, we always use $\xi$ to denote 
a fixed primitive element of $\Ff_q$.
For $a,b,c\in\Ff_q$, we define
\begin{equation*}
{\bf d}(a,b,c):=
\begin{pmatrix}
a & 0 & 0 \\
0 & b & 0 \\
0 & 0 & c\\
\end{pmatrix},\;
{\bf ad}(a,b,c):=
\begin{pmatrix}
0& 0 & a\\
0 & b & 0\\
c & 0 & 0\\
\end{pmatrix}.
\end{equation*}
For the convenience of discussion, we adopt the following
special representatives of $G$ from \cite{swx}:
\begin{equation}\label{groupg}
H\cup {\bf d}(1,\xi^{-1},\xi^{-2})\cdot H.
\end{equation}
Moreover, the following holds.
\begin{Lemma}\label{transitive}\cite{dye}
The group $G$ acts transitively on both $I$
$($respectively, $Pa)$ and $E$ $($respectively, $Se)$.
\end{Lemma}
\noindent{The following result was proved in \cite{swx} and
will be used frequently.}

\begin{Lemma}\label{meet}\cite[Lemma 2.9]{swx}
Let $\p$ be a point not on ${\mathcal O}$, $\ell$ a non-tangent
line, and $\p\in\ell$. Using the above notation, we have the
following.

\begin{enumerate}
\renewcommand{\labelenumi}{(\roman{enumi})}

\item If $\p\in I$ and $\ell\in$Pa, then $\p^\perp\cap \ell \in E$
if $q\equiv 1 \pmod 4$, and $\p^\perp \cap \ell \in I$ if $q\equiv 3
\pmod 4$.

\item If $\p\in I$ and $\ell \in$Se, then $\p^\perp\cap \ell \in I$
if $q\equiv 1 \pmod 4$, and $\p^\perp\cap \ell \in E$ if $q\equiv 3
\pmod 4$.

\item If $\p\in E$ and $\ell \in$Pa, then $\p^\perp\cap \ell \in I$
if $q\equiv 1 \pmod 4$, and $\p^\perp\cap \ell \in E$ if $q\equiv 3
\pmod 4$.

\item If $\p\in E$ and $\ell \in$Se, then $\p^\perp\cap \ell\in E$
if $q\equiv 1 \pmod 4$, and $\p^\perp\cap\ell\in I$ if $q \equiv 3
\pmod 4$.

\end{enumerate}
\end{Lemma}


Next we define $\Sq_q-1:=\{s-1\mid s\in \Sq_q\}$
and $\Nsq_q-1:=\{s-1\mid s\in \Nsq_q\}$.  
\begin{Lemma}\label{cs}\cite{store}
Using the above notation,
\begin{itemize}
\item[(i)] if $q\equiv 1\pmod 4$, then 
$|(\Sq_q-1)\cap\Sq_q|=\frac{q-5}{4}$
and $$|(\Sq_q-1)\cap \Nsq_q|=|(\Nsq_q-1)\cap\Sq_q|=|(\Nsq_q-1)\cap\Nsq_q|=\frac{q-1}{4};$$
\item[(ii)] if $q\equiv 3\pmod 4$, 
then $|(\Nsq_q-1)\cap\Sq_q|=\frac{q+1}{4}$ 
and $$|(\Sq_q-1)\cap\Sq_q|=|(\Sq_q-1)\cap\Nsq_q|=|(\Nsq_q-1)\cap\Nsq_q|
=\frac{q-3}{4}.$$
\end{itemize}
\end{Lemma}

\begin{Definition}
Let $\p$ be a point not on $\mathcal{O}$ and $\ell$ a line. 
We define $E_{\ell}$ and $I_{\ell}$ to be the set of external 
points and the set of internal points on $\ell$, respectively,
$Pa_{\p}$ and $Se_{\p}$ the set of passant lines and the set 
of secant lines through $\p$, respectively, and $T_{\p}$ the 
set of tangent lines through $\p$. Also, $N(\p)$ is defined 
to be the set of internal points on the passant lines through 
$\p$ including or excluding $\p$ accordingly as 
$q\equiv 3\pmod 4$ or $q\equiv 1\pmod 4$. 
\end{Definition}

\begin{Remark}\label{bsize}
Using the above notation, for $\p\in I$, we have
$|E_{\p^\perp}| = |Se_{\p}| = \frac{q+1}{2}$;
$|I_{\p^\perp}| = |Pa_{\p}| = \frac{q+1}{2}$; and 
$|N(\p)| = \frac{q^2-1}{4}$ or $\frac{q^2+3}{4}$
accordingly as $q\equiv 1\pmod 4$ or $q\equiv  3\pmod 4$.
\end{Remark}
Let $\p\in I$, $\ell\in Pa$, $g\in G$, and $W\le G$.
Using standard notations from permutation group theory,
we have 
$I_{\ell}^g=I_{\ell^g}$, $Pa_{\p}^g=Pa_{\p^g}$;
$E_{\ell}^g=E_{\ell^g}$, $Se_{\p}^g=Se_{\p^g}$;
$H_{\p}^g=H_{\p^g}$;
$N(\p)^g = N(\p^g)$, $(W^g)_{\p^g} = W_{\p}^g$. 
We will use these results later without further 
reference. Also, the definition of $G$ yields that 
$(\p^\perp)^g=(\p^g)^\perp$, where $\perp$ is 
the above defined polarity of $\PG(2,q)$.




\begin{Proposition}\label{Ktransitive}
Let $\p\in I$ and set $K:=G_{\p}$. Then $K$ is transitive on  
$I_{\p^\perp}$, $E_{\p^\perp}$, $Pa_{\p}$, and $Se_{\p}$, 
respectively.
\end{Proposition}
{\Proof} Witt's theorem \cite{bh} implies that $K$ acts
transitively on isometry classes of the form $Q$
on the points of $\p^\perp$. Note that $K=G_{\p^\perp}$
by the definition of $G$. Dually, we must have that
$K$ is transitive on both $Pa_{\p}$ and $Se_{\p}$.
\QED
When $\p=(1,0,-\xi)$, using (\ref{lemB}), (\ref{grouph}),
and (\ref{groupg}), we obtain that $K:= G_\p=$
\begin{equation}\label{stabP}
\begin{array}{lllllll}
\left.\left\{\left(\begin{array}{ccc} d^2 & cd\xi & c^2 \xi^2\\ 
2cd & d^2+c^2\xi & 2dc\xi\\ c^2 & dc & d^2 \end{array}\right)\right| 
d,c \in \Ff_q, d^2-c^2\xi =1\right\}\\
\bigcup\left.\left\{\left(\begin{array}{ccc} d^2 & -cd\xi & c^2 \xi^2\\ 
2cd& -d^2-c^2\xi& 2dc\xi\\ c^2 & -dc & d^2 \end{array}\right)\right| 
d,c \in \Ff_q, -d^2+c^2\xi =1\right\}\\
\bigcup\left.\left\{\left(\begin{array}{ccc} d^2 & cd & c^2\\ 
2cd\xi^{-1} & d^2+c^2\xi^{-1}& 2dc\\ c^2\xi^{-2} & dc\xi^{-1} & d^2 
\end{array}\right)\right| d,c \in \Ff_q, d^2\xi-c^2 =1\right\}\\
\bigcup\left.\left\{\left(\begin{array}{ccc} d^2 & -cd & c^2 \\ 
2cd\xi^{-1} & -d^2-c^2\xi^{-1}& 2dc\\ c^2\xi^{-2} & -dc\xi^{-1} & d^2 
\end{array}\right)\right| d,c \in \Ff_q, -d^2\xi+c^2 =1\right\}.\\
\end{array}
\end{equation}
Recall that $N(\p)$ for $\p\in I$ is the set of internal points on the
passant lines through $\p$, where $\p$ is included or not accordingly
as $q\equiv 3\pmod 4$ or $q\equiv 1\pmod 4$.

\begin{Theorem}\label{intersection1}
Let $\p\in I$ and $\ell\in Pa$. Then $|N(\p) \cap I_{\ell}|\equiv 0\pmod 2$.
\end{Theorem}
{\Proof} If $\p\in \ell$, it is clear that
\begin{equation*}
|N(\p)\cap I_{\ell}|=
\begin{cases}
\frac{q-1}{2}, &\text{if}\; q\equiv 1\pmod 4,\\
\frac{q+1}{2}, & \text{if}\;q\equiv 3\pmod 4,
\end{cases}
\end{equation*}
which is even. Therefore, $|N(\p)\cap I_{\ell}|\equiv 0 \pmod 2$ for
this case.

If $\ell=\p^\perp$, by Lemma~\ref{meet}(i), we have
\begin{equation*}
|N(\p)\cap I_{\ell}|=
\begin{cases}
0, &\text{if}\;q\equiv 1\pmod 4,\\
\frac{q+1}{2}, &\text{if}\;q\equiv 3\pmod 4,
\end{cases}
\end{equation*}
which is even. Hence, $|N(\p)\cap I_{\ell}|\equiv 0 \pmod 2$ 
for this case.

Now we assume that neither $\ell=\p^\perp$ nor $\p\in\ell$. As
$G$ is transitive on $Pa$ and preserves incidence, we may take 
$\ell=\p_1^\perp=[1,0,-\xi^{-1}]$, where $\p_1=(1,0,-\xi)\in I$.
Since $\p$ is either on a passant line through $\p_1$ or on a 
secant line through $\p_1$, the rest it to show that 
$|N(\p)\cap I_{\ell}|$ is even for any $\p$ on a line through
$\p_1$ with $\p\notin \ell$ and $\p\not=\p_1$. 

{\bf Case I.} $\p$ is a point on a secant line through 
$\p_1$ and $\p\notin \ell$.

Since $K=G_{\p_1}$ acts transitively on $Se_{\p_1}$ 
by Proposition~\ref{Ktransitive}, it is enough to establish 
that $|N(\p)\cap I_{\ell}|$ is even for an arbitrary internal 
point on a {\it special} secant line, $\ell_1$ say, through 
$\p_1$. To this end, we may take $\ell_1=[0,1,0]$.
It is clear that
$$I_{\ell_1}=\{(1,0,-\xi^j)\mid 0\le j\le q-1,\;\text{$j$ odd}\}$$
and
$$I_{\ell}=\{(1,s,\xi)\mid s \in \Ff_q,\; s^2-\xi\in\Nsq_q \}.$$
Hence, if $\p=(1,0,-\xi^j)\in I_{\ell_1}$ then 
$$D_j=\left.\left\{\left[1,-\frac{\xi^{1-j}+1}{s},\frac{1}{\xi^j}\right]
\right|s\in\Ff_q^*,\;s^2-\xi\in\Nsq_q\right\}\cup\{[0,1,0]\}$$
consists of the lines through both $\p$ and the points on $\ell$. 
Note that the number of passant lines in $D_j$ is determined 
by the number of $s$ satisfying both
\begin{equation}\label{eq1}
\frac{1}{s^2}(\xi^{1-j}+1)^2-\frac{4}{\xi^j}
\in \Nsq_q
\end{equation}
and
\begin{equation}\label{eq2}
s^2-\xi\in \Nsq_q.
\end{equation}
Since, $s\not=0$ and whenver $s$ satisfies both (\ref{eq1}) 
and (\ref{eq2}), so does $-s$, we see that $|N(\p)\cap I_{\ell}|$ 
must be even in this case.

{\bf Case II.} $\p$ is an internal point on a passant 
line through $\p_1$ and $\p\notin\ell$.

By Lemma~\ref{meet}, we may assume that 
$\p\in\p_3^\perp$, where $\p_3=(1,x,\xi)\in I_{\ell}$ 
with $x\in\Ff_q^*$ and $x^2-\xi\in\Nsq_q$.
Here $\p_3^\perp=[1,-\frac{2x}{\xi},\frac{1}{\xi}]$ is 
a passant line through $\p_1$. Let $K=G_{\p_1}$. Using 
(\ref{stabP}), we can obtain $L:=K_{\p_3}$ as follows: 
if $q\equiv 1\pmod 4$,
\begin{equation}\label{Sub_Stab_I}
\begin{array}{llll}
L & = &
\left\{
\left(
\begin{array}{rrr}
\frac{\xi}{x^2-\xi} & \frac{x\xi}{x^2-\xi} & \frac{x^2\xi}{x^2-\xi}\\
-\frac{2x}{x^2-\xi} & -\frac{x^2+\xi}{x^2-\xi} & -\frac{2x\xi}{x^2-\xi}\\
\frac{x^2}{(x^2-\xi)\xi} & \frac{x}{x^2-\xi} & \frac{\xi}{x^2-\xi}
\end{array}
\right),
\left(
\begin{array}{rrr}
0 & 0 & -1\\
0 & -\frac{1}{\xi} & 0 \\
-\frac{1}{\xi^2} & 0 & 0\\
\end{array}
\right)
\right\}\\
{} &\bigcup&
\left\{
\left(
\begin{array}{rrr}
-\frac{x^2}{(x^2-\xi)\xi} & -\frac{x}{x^2-\xi} & -\frac{\xi}{x^2-\xi}\\
\frac{2x}{(x^2-\xi)\xi} & \frac{x^2+\xi}{(x^2-\xi)\xi} & \frac{2x}{x^2-\xi}\\
-\frac{1}{(x^2-\xi)\xi} & -\frac{x}{(x^2-\xi)\xi} & -\frac{x^2}{(x^2-\xi)\xi}\\
\end{array}
\right),
\left(
\begin{array}{rrr}
1 & 0 & 0\\
0 & 1 & 0 \\
0 & 0 &  1\\
\end{array}
\right)
\right\};
\end{array}
\end{equation}
if $q\equiv 3\pmod 4$,
\begin{equation}\label{Sub_Stab_I_1}
\begin{array}{llll}
L &
= &\left\{
\left(
\begin{array}{rrr}
\frac{\xi}{x^2-\xi} & \frac{x\xi}{x^2-\xi} & \frac{x^2\xi}{x^2-\xi}\\
-\frac{2x}{x^2-\xi} & -\frac{x^2+\xi}{x^2-\xi} & -\frac{2x\xi}{x^2-\xi}\\
\frac{x^2}{(x^2-\xi)\xi} & \frac{x}{x^2-\xi} & \frac{\xi}{x^2-\xi}\\
\end{array}
\right), 
\left(
\begin{array}{rrr}
0 & 0 & -\xi\\
0 & -1 & 0 \\
-\frac{1}{\xi} & 0 & 0\\
\end{array}
\right)
\right\}\\
{} &\cup &\left\{
\left(\begin{array}{rrr}
-\frac{x^2}{x^2-\xi} & -\frac{x\xi}{x^2-\xi} & -\frac{\xi^2}{x^2-\xi}\\
\frac{2x}{x^2-\xi} & \frac{x^2+\xi}{x^2-\xi} & \frac{2x\xi}{x^2-\xi}\\
-\frac{1}{x^2-\xi} & -\frac{x}{x^2-\xi} & -\frac{x^2}{x^2-\xi}
\end{array}
\right),
\left(
\begin{array}{rrr}
1 & 0 & 0 \\
0 & 1 & 0 \\
0 & 0 & 1 \\
\end{array}
\right)
\right\}.
\end{array}
\end{equation}

Let $(1,y,\xi)$ be a point on ${\ell}$. Using
(\ref{Sub_Stab_I}) and (\ref{Sub_Stab_I_1}), 
we have that $L$ fixes $(1,y,\xi)$ if and only if 
$$xy^2-(x^2+\xi)y+x\xi=0;$$
that is, $y=x$ or $y=\frac{\xi}{x}$. Consequently, 
$\p_3=(1,x,\xi)$ and $\ell\cap \p_3^\perp=(1,\frac{\xi}{x}, \xi)$
are the only points of the form $(1,s,t)$ on $\ell$
fixed by $L$.
Since $\p\in\p_3^\perp$, $\p\not=\p_1$ and 
$\p\not=\p_3^\perp\cap \ell$,
$\p=(1,\frac{\xi+n}{2x}, n)$ for some $n\not=\xi$. 
Now if we denote by $\mathbf{V}$ the set of passant 
lines through 
$\p$ that meet $\ell$ in an internal point, then it 
is clear that $|\mathbf{V}|=|N(\p)\cap I_{\ell}|$. 
Direct computations give us that 
$L_{\p}\cong \mathbb{Z}_2$. Since $\p_3$ and $\p$ 
are both fixed by $L_{\p}$, it follows that both 
$\ell_{\p_3,\p}$ and $\p_3^\perp$ are fixed by $L_{\p}$. 
Note that when $q\equiv 3\pmod 4$, both $\p_3^\perp$ 
and $\ell_{\p_3,\p}$ are in $\mathbf{V}$; and when 
$q\equiv 1\pmod 4$, neither $\ell_{\p_3,\p}$ nor 
$\p_3^\perp$ is in $\mathbf{V}$. If there were another 
line $\ell^{'}$ through $\p$ which is distinct from 
both $\p_3^\perp$ and $\ell_{\p_3,\p}$ and which is 
also fixed by $L_{\p}$, then $L_{\p}$ would fix at 
least three points on $\ell=\p^\perp$, namely, 
$\ell^{'}\cap \ell$, $\p_3^\perp\cap \ell$, and $\p_3$. 
Since no further point of the form $(1,s,t)$ except 
for $\p_3$ and $\ell\cap \p_3^\perp$ can be fixed by 
$L$ due to the above discussion, we must have 
$\ell^{'}\cap \ell=(0,1,0)\in E_{\ell}$.
So $\ell^{'}\notin\mathbf{V}$. Using the fact that 
$L_\p$ preserves incidence, we conclude that when 
$q\equiv 1\pmod 4$, $L_\p$ has $\frac{|\mathbf{V}|}{2}$ 
orbits of length $2$ on $\mathbf{V}$; and when 
$q\equiv 3\pmod 4$, $L_{\p}$ has two orbits of 
length $1$, namely, $\{\p_3^\perp\}$ and 
$\{\ell_{\p_3,\p}\}$, and  $\frac{|\mathbf{V}|-2}{2}$ 
orbits of length $2$ on $\mathbf{V}$. Either forces 
$|\mathbf{V}|$ to be even. Therefore, 
$|N(\p)\cap I_{\ell}|$ is even.
\QED
Recall that $\A$ is the incidence matrix of $Pa$ and $I$ 
whose columns are indexed by the internal points $\p_1$, 
$\p_2$,..., $\p_N$ and whose rows are indexed by the passant 
lines $\p_1^\perp$, $\p_2^\perp$,..., $\p_N^\perp$; and 
$\A$ is symmetric. For the convenience of discussion, 
for $\p\in I$, we define
$$
\widehat{N(\p)}=
\begin{cases}
N(\p)\cup\{\p\}, & \textrm{if}\;q\equiv 1\pmod 4,\\
N(\p)\setminus\{\p\}, & \textrm{if}\;q\equiv 3\pmod 4.\\
\end{cases}
$$
That is, $\widehat{N(\p)}$ is the set of the internal
points on the passant lines through $\p$ including
$\p$. It is clear that for $\p\notin\ell$,
$|N(\p)\cap I_\ell| = |\widehat{N(\p)}\cap I_\ell|$.
\begin{Lemma}
Using the above notation, we have $\A^3\equiv \A \pmod 2$, 
where the congruence means entry-wise congruence.
\end{Lemma}
{\Proof} Since the $(i,j)$-entry of $\A^2=\A^\top\A$ is 
the standard dot product of the $i$-th row of $\A^\top$ 
and $j$-th column of $\A$, we have
\begin{equation*}
(\A^2)_{i,j}=(\A^\top\A)_{i,j}=
\begin{cases}
\frac{q+1}{2}, & \text{if}\;i=j,\\
1, & \text{if}\;\ell_{\p_i, \p_j}\in Pa, \\
0, & \text{otherwise}.
\end{cases}
\end{equation*}
Therefore, the $i$-th row of $\A^2\pmod 2$ indexed by 
$\p_i$ can be viewed as the characteristic row vector 
of $\widehat{N(\p_i)}$.

If $\p_i\in \p_j^\perp$, then 
$(\A^3)_{i,j}=((\A^2)\A^\top)_{i,j}=q$ since 
$(\A^2)_{i,i}=\frac{q+1}{2}$ and there are 
$\frac{q-1}{2}$ internal points other than $\p_i$ on 
$\p_j^\perp$ that are connected with $\p_i$ by the 
passant line $\p_j^\perp$. If $\p_i\not\in \p_j^\perp$, 
then 
$(\A^3)_{i,j}=((\A^\top\A)\A^\top)_{i,j}\equiv
|\widehat{N(\p_i)}\cap I_{\p_j^\perp}|=
|N(\p_i)\cap I_{\p_j^\perp}|\equiv 0\pmod 2$
by Theorem~\ref{intersection1}. Consequently,
\begin{equation*}
(\A^3)_{i,j}\equiv
\begin{cases}
1\pmod 2, & \text{if}\; \p_i\in \p_j^\perp,\\
0 \pmod 2, & \text{if}\;\p_i\notin\p_j^\perp.
\end{cases}
\end{equation*}
The lemma follows immediately. \QED

\section{The Conjugacy Classes and Intersection Parity}
In this section, we present detailed information about 
the conjugacy classes of $H$ and study their 
intersections with some special subsets of $H$.
\subsection{Conjugacy classes}
The conjugacy classes of $H$ can be read off in 
terms of the map $T=\textrm{tr}(g)+1$, where 
$\textrm{tr}(g)$ is the trace of $g$.


\begin{Lemma}\cite[Lemma 3.2]{swx}\label{classes}
The conjugacy classes of $H$ are given as follows.
\begin{itemize}
\item[(i)] $D=\{{\bf d}(1,1,1)\}$;
\item[(ii)] $F^{+}$ and $F^{-}$, where $F^{+}\cup 
F^{-} = \{g\in H\mid T(g) = 4,\;g\not={\bf d}(1,1,1)\}$;
\item[(iii)] $[\theta_i] =\{g\in H \mid T(g) = 
\theta_i\}$, $1\le i \le \frac{q-5}{4}$ if 
$q\equiv 1\pmod 4$, or $1\le i\le \frac{q-3}{4}$ 
if $q\equiv 3 \pmod 4$, where $\theta_i\in \Sq_q$, 
$\theta_i\not= 4$, and $\theta_i-4\in\Sq_q$;
\item[(iv)] $[0]=\{g\in H\mid T(g) =0\}$;
\item[(v)] $[\pi_k]=\{g\in H\mid T(g) = \pi_k\}$, 
$1\le k\le \frac{q-1}{4}$ if $q\equiv 1\pmod 4$, 
or $1\le k\le \frac{q-3}{4}$ if $q\equiv 3 \pmod 4$, 
where $\pi_i\in\Sq_q$, $\pi_k\not=4$, and 
$\pi_k-4\in\Nsq_q$.
\end{itemize}
\end{Lemma}
\begin{Remark}\label{order}
The set $F^+\cup F^-$
forms one conjugacy class of $G$, and splits into two 
equal-sized classes $F^{+}$ and $F^{-}$ of $H$. For our 
purpose, we denote $F^+\cup F^-$ by $[4]$. Also, each 
of $D$, $[\theta_i]$, $[0]$, and $[\pi_k]$ forms a single 
conjugacy class of $G$. The class $[0]$ consists of all 
the elements of order $2$ in $H$.
\end{Remark}
In the following, for convenience, we frequently use $C$ 
to denote any one of $D$, $[0]$, $[4]$, $[\theta_i]$, or
$[\pi_k]$. That is,
\begin{equation}\label{defC}
C=D, [0], [4], [\theta_i], \;{\rm or}\;[\pi_k].
\end{equation}

\subsection{Intersection properties}

\begin{Definition}
Let $\p,\q \in I$, $W\subseteq I$, and $\ell \in Pa$. We 
define $\mathcal{H}_{\p,\q}=\{h\in H\mid (\p^\perp)^h\in Pa_{\q}\}$,
$\mathcal{S}_{\p,\ell}=\{h\in H\mid (\p^\perp)^h=\ell\}$, 
and $\mathcal{U}_{\p, W}=\{h\in H\mid \p^h\in W\}$.
That is, $\mathcal{H}_{\p,\q}$ consists of all the elements 
of $H$ that map the passant line $\p^\perp$ to a passant 
line through $\q$, $\mathcal{S}_{\p,\ell}$ is the set of 
elements of $H$ that map $\p^\perp$ to the passant line 
$\ell$, and $\mathcal{U}_{\p,W}$ is the set of elements 
of $H$ that map $\p$ to a point in $W$.
\end{Definition}
Using the above 
notation, we have that 
$\mathcal{H}_{\p,\q}^g=\mathcal{H}_{\p^g,\q^g}$,
$\mathcal{S}_{\p,\ell}^g=\mathcal{S}_{\p^g,\ell^g}$, 
and $\mathcal{U}_{\p,W}^g=\mathcal{U}_{\p^g,W^g}$, 
where $\mathcal{H}_{\p,\q}^g=\{g^{-1} h g\mid h\in \mathcal{H}_{\p,\q}\}$,
$\mathcal{S}_{\p,\ell}^g=\{g^{-1} h g\mid h\in \mathcal{S}_{\p,\q}\}$,
and $\mathcal{U}_{\p,W}^g=\{h^g\mid h\in\mathcal{U}_{\p,W}\}$. 
Moreover, it is true that 
$(C\cap \mathcal{H}_{\p,\q})^g=C\cap \mathcal{H}_{\p^g, \q^g}$
and $(C\cap \mathcal{U}_{\p, W})^g=C\cap \mathcal{U}_{\p^g, W^g}$. 
In the following discussion, we will use these results 
without further reference.

Next we compute the size of the intersection of each conjugacy 
class of $H$ with $K$ which is a stabilizer of an internal 
point in $H$.
\begin{Corollary}\label{Kintersection}
Let $\p \in I$ and $K=H_{\p}$. Then we have
\begin{itemize}
\item[(i)] $|K\cap D|=1$;
\item[(ii)] $|K\cap [4]| = 0$;
\item[(iii)] $|K\cap [\pi_k]| = 2$;
\item[(iv)] $|K\cap [\theta_i]| = 0$;
\item[(v)] $|K\cap [0]|=\frac{q+1}{2}$ or $\frac{q-1}{2}$ 
accordingly as $q\equiv 1\pmod 4$ or $q\equiv 3\pmod 4$.
\end{itemize}
\end{Corollary}
{\Proof} Let $\q=(1,0,-\xi)$ and $K_1=H_{\q}$. As $H$ 
is transitive on $I$, $\q^g=\p$ for some $g\in H$. 
Moreover, $$|K\cap C|=|K_1^g\cap C|=|(K_1\cap C)^g|.$$
Therefore, to prove the corollary, it is 
sufficient to assume that $\p=\q$.

It is obvious that $|D\cap K|=1$. Let $g\in K\cap C$. 
Using (\ref{stabP}), we have that the quadruples 
$(a,b,c,d)$ determining $g$ satisfy the following 
equations
\begin{equation}\label{linear1}
\begin{array}{ccccc}
bd-ac\xi& = & 0\\
b^2-a^2\xi& = & -\xi(d^2-c^2\xi)\\
ad-bc& = & 1\\
a+d & = & s,
\end{array}
\end{equation}
where $s^2=0$, $4$, $\pi_k$, $\theta_i$. The equations 
in (\ref{linear1}) yield (1) $a=d=\frac{s}{2}$, 
$c^2=\frac{s^2-4}{4\xi}$, $b^2=\frac{(s^2-4)\xi}{4}$ and 
(2) $a=-d$, $s=0$, $c^2\xi-1=a^2$. From Case (1), we have 
$|K\cap [\pi_k]|=2$ for each $[\pi_k]$ and $|K\cap C|=0$ 
for $C=[\theta_i]$, $[4]$; moreover, if $q\equiv 3\pmod 4$, 
we obtain one group element 
{\bf ad}$(-\xi, -1,\xi^{-1})\in K\cap [0]$
in Case (1). Since the number of $t\in\Nsq_q$ satisfying 
$t-1\in \Sq_q$ is $\frac{q-1}{4}$ or $\frac{q-3}{4}$ 
accordingly as $q\equiv 1\pmod 4$ or $q\equiv 3\pmod 4$ 
by Lemma~\ref{bsize}, the number of $c\in\Ff_q^*$ 
satisfying $c^2\xi-1\in\Sq_q$ is $2|(\Nsq_q-1)\cap\Sq_q|$ 
which is $\frac{q-1}{2}$ or $\frac{q-3}{2}$ accordingly 
as $q\equiv 1\pmod 4$ or $q\equiv 3\pmod 4$. 
When $q\equiv 1\pmod 4$, $c=0$ also satisfies 
$c^2\xi-1\in\Sq_q$. Therefore, Case (1) and Case (2) 
give rise to $\frac{q+1}{2}$ or $\frac{q-1}{2}$ different 
group elements in $K\cap [0]$ depending on $q$. The 
corollary now is proved.
\QED

In the following lemmas, we investigate the parity of 
$|\mathcal{H}_{\p,\q}\cap C|$ for $C\not=[0]$ and 
$\p,\q\in I$. Recall that $\ell_{\p,\q}$ is the line 
through $\p$ and $\q$.
\begin{Lemma}\label{m1}
Let $\p,\q\in I$. Suppose that $C=D$,
$[4]$, $[\pi_k]$ $(1\le k\le \frac{q-1}{4})$, 
or $[\theta_i]$ $(1\le i\le \frac{q-5}{4})$.

First assume that $q\equiv 1\pmod 4$.
\begin{itemize}
\item[(i)] If $\ell_{\p,\q}\in Pa_{\p}$, then 
$|\mathcal{H}_{\p,\q}\cap C|$ is always even.
\item[(ii)] If $\ell_{\p,\q}\in Se_{\p}$, 
$\q\notin \p^\perp$, and $|\mathcal{H}_{\p,\q}\cap C|$ 
is odd, then $C=[\theta_{i_1}]$ or $[\theta_{i_2}]$.
\item[(iii)] If $\q\in\ell_{\p,\q}\cap \p^\perp$ and
$|\mathcal{H}_{\p,\q}\cap C|$ is odd, then $C=D$. 
\end{itemize}

Now assume that $q\equiv 3\pmod 4$.
\begin{itemize}
\item[(iv)] If $\ell_{\p,\q}\in Se_{\p}$, then 
$|\mathcal{H}_{\p,\q}\cap C|$ is always even.
\item[(v)] If $\ell_{\p,\q}\in Pa_{\p}$, $\q\notin \p^\perp$, 
and $|\mathcal{H}_{\p,\q}\cap C|$ is odd, then $C=[\pi_{i_1}]$ 
or $[\pi_{i_2}]$.
\item[(vi)] If $\q\in\ell_{\p,\q}\cap \p^\perp$ and
$|\mathcal{H}_{\p,\q}\cap C|$ is odd, then $C=D$. 
\end{itemize}
\end{Lemma}

{\Proof} We only provide the detailed proof for the case 
when $q\equiv 1\pmod 4$. Since $G$ acts transitively on 
$I$ and preserves incidence, without loss of generality, 
we may assume that $\p=(1,0,-\xi)$ and let $K=G_{\p}$.

Since $K$ is transitive on both $Pa_{\p}$ and 
$Se_{\p}$ by Proposition~\ref{Ktransitive} and 
$|\mathcal{H}_{\p,\q}\cap C|$ = 
$|(\mathcal{H}_{\p,\q}\cap C)^g|$=$|\mathcal{H}_{\p,\q^g}\cap C|$, 
we may assume that $\q$ is on either $\ell_1$ or 
$\ell_2$, where
$\ell_1=[1,0,\xi^{-1}]\in Pa_{\p}$ and $\ell_2=[0,1,0]\in Se_{\p}$.

{\bf Case I.} $\q\in \ell_1$.

In this case, $\q=(1,x,-\xi)$ for some $x\in\Ff_q^*$ 
and $x^2+\xi\in\Nsq_q$, and 
$$Pa_{\q}=\{[1,s,(1+sx)\xi^{-1}]\mid s\in \Ff_q, s^2-4(1+sx)\xi^{-1}\in\Nsq_q\}.$$
Using (\ref{stabP}), we obtain that
$$K_{\q}=\{{\bf d}(1,1,1), {\bf ad}(1,-\xi^{-1},\xi^{-2})\}.$$
It is obvious that ${\bf d}(1,1,1)$ fixes each line 
in $Pa_{\q}$. From
$${\bf ad}(1,-\xi^{-1},\xi^{-2})^{-1}(1,s,(1+sx)\xi^{-1})^\top=((1+sx)\xi, -s\xi, 1)^\top,$$
it follows that a line of the form $[1,s,(1+sx)\xi^{-1}]$ 
is fixed by $K_{\q}$ if and only if $s=0$ or $s=-2x^{-1}$. 
Further, since $[1,-2x^{-1}, -\xi^{-1}]$ is a secant line, 
we obtain that $K_{\q}$ on $Pa_{\q}$ has one orbit of 
length $1$, i.e. $\{\ell_1=[1,0,\xi^{-1}]\}$, and all 
other orbits, whose representatives are $\mathcal{R}_1$, 
have length $2$. From
\begin{equation*}
|\mathcal{H}_{\p,\q}\cap C|=|\mathcal{S}_{\p,\ell_1}\cap C|
+
2\displaystyle\sum_{\ell\in\mathcal{R}_1}|\mathcal{S}_{\p,\ell}\cap C|,
\end{equation*}
it follows that the parity of $|\mathcal{H}_{\p,\q}\cap C|$ 
is determined by that of $|\mathcal{S}_{\p,\ell_1}\cap C|$. 
Here we used the fact that $|\mathcal{S}_{\p,\ell}\cap C|
=|\mathcal{S}_{\p,\ell^{'}}\cap C|$ if $\{\ell, \ell^{'}\}$ 
is an orbit of $K_{\p}$ on $Pa_{\q}$. Meanwhile, it is clear 
that $|\mathcal{S}_{\p,\ell_1}\cap D|=0$.

Note that the quadruples $(a,b,c,d)$ that determine group 
elements in $\mathcal{S}_{\p,\ell_1}\cap C$ are the solutions 
to the following equations
\begin{equation}\label{s1}
\begin{array}{ccccc}
-2cd+2ab\xi^{-1} & = & 0\\
c^2-a^2\xi^{-1} & = & (d^2-b^2\xi^{-1})\xi^{-1}\\
a+d & = & s\\
ad-bc & = & 1,\\
\end{array}
\end{equation}
where $s^2=4$, $\pi_k$, $\theta_i$, and that if one of
$b$ and $c$ is zero, so is the other. If $b=c=0$ and
$2\in\Sq_q$ then the above (\ref{s1}) give $4$ group 
elements in $[2]$ and $0$ element is any other classes.
 If neither $b$ nor $c$ is zero, then the first 
two equations in (\ref{s1}) yield $b=\pm\sqrt{-1}\xi c$. 
Combining with the last two equations in (\ref{s1}), 
we obtain $0$, $4$ or $8$ quadruples $(a,b,c,d)$ 
satisfying the above equations, among which both 
$(a,b,c,d)$ and $(-a,-b,-c,-d)$ appear at the same 
time. Since $(a,b,c,d)$ and $(-a,-b,-c,-d)$ give 
rise to the same group element, we conclude that 
$|\mathcal{S}_{\p,\ell_1}\cap C|$ is $0$, $2$, 
or $4$.

{\bf Case II.} $\q\in \ell_2$, $\q\notin \p^\perp$, 
and $\q\not=\p$.

In this case, $\q=(1,0,-y)$ for some $y\in\Nsq_q$ 
and $y\not=\pm\xi$. Using (\ref{stabP}), we obtain 
that
$$K_{\q}=\{{\bf d}(1,1,1), {\bf d}(-1,1,-1)\}.$$
Moreover, $K_{\q}$ on 
$Pa_{\q}=\{[1,s,y^{-1}]\mid s\in\Ff_q, s^2-4y^{-1}\in\Nsq_q\}$ 
has one orbit of length $1$, that is, 
$\{\ell_4=[1,0,y^{-1}]\}$, and all other orbits are 
of length $2$. Similar arguments as above show 
that the parity of $|\mathcal{H}_{\p, \q}\cap C|$ 
is the same as that of $|\mathcal{S}_{\p,\ell_4}\cap C|$. 
So the rest is to find the parity of 
$|\mathcal{S}_{\p,\ell_4}\cap C|$. The group elements 
in $\mathcal{S}_{\p,\ell_4}\cap C$ are determined
by the quadruples $(a,b,c,d)$ satisfying the following 
equations:
\begin{equation}\label{text1}
\begin{array}{cccc} 
-2cd + 2ab\xi^{-1} & = & 0\\
c^2-a^2\xi^{-1} & = & (d^2-b^2\xi^{-1})y^{-1}\\
a+d & = & s\\
ad-bc & = & 1.\\
\end{array}
\end{equation}
Note that if one of $b$ and $c$ is zero, so is the
other.
If neither $b$ nor $c$ is zero, then the first two 
equations in (\ref{text1}) yield $b=\pm\sqrt{-\xi y}c$ 
and $a=\pm\sqrt{-\xi y^{-1}} d$. Combining with the 
last two, 
the above quadruples $(a,b,c,d)$ yield 
$0$, $2$, or $4$ group elements in $[s^2]$. If $b=c=0$, then
$ad=1$, $d^2=\pm\sqrt{-y\xi^{-1}}$ and $a^2=\pm\sqrt{-\xi y^{-1}}$;
and so 
$$s^2=\sqrt{-\xi y^{-1}}+\sqrt{- y \xi^{-1}}+2
\;\;\;\textrm{or}\;\;\;s^2=
-\sqrt{-\xi y^{-1}}-\sqrt{-y\xi^{-1}}+2.$$
Since $(\sqrt{-\xi y^{-1}}+\sqrt{-y\xi^{-1}}+2)
(-\sqrt{-\xi y^{-1}}-\sqrt{-y \xi^{-1}}+2)=
(\sqrt{\xi y^{-1}}+\sqrt{y\xi^{-1}})^2$, 
the above quadruples $(a,b,c,d)$ yield
$0$ or $1$ group elements in two classes $[\theta_{i_1}]$ 
and $[\theta_{i_2}]$ where 
$\theta_{i_1}=\sqrt{-\xi y^{-1}}+\sqrt{- y \xi^{-1}}+2$ and
$\theta_{i_2}=-\sqrt{-\xi y^{-1}}-\sqrt{-y\xi^{-1}}+2$. 
The above analysis shows that if $|\mathcal{H}_{\p,\q}\cap C|$ 
is odd then $C=[\theta_{i_1}]$ or $[\theta_{i_2}]$ in this case.

{\bf Case III.} $\q=\ell_2\cap \p^\perp$.

In this case, $\q=(1,0,\xi)$ and the set of passant lines
through $\q$ is $$Pa_{\q}=\{[1,u,-\xi^{-1}]\mid u \in \Ff_q, 
u^2+\xi\in\Nsq_q\}.$$ Using (\ref{stabP}), we obtain that
$$K_\q=\{{\bf d}(1,1,1), {\bf d}(-1,1,-1),{\bf ad}(-1,-\xi^{-1},
-\xi^{-2}), {\bf ad}(1,-\xi^{-1},\xi^{-2})\}.$$
Therefore, among the orbits of $K_\q$ on $Pa_\q$, 
$\{[1,0,-\xi^{-1}]\}$ is the only one
of length $1$ and all others are of length $2$. Hence,
the parity of $|\mathcal{H}_{\p,\q}\cap C|$ is the same as
that of $|\mathcal{S}_{\p,\p}\cap C|$ which is the same as 
that of $|K\cap C|$; by Corollary~\ref{Kintersection}, it 
follows that $|K\cap C|$ is odd if and only if $C=D$.

\QED
For $\q\in I$, we denote by $\overline{N(\q)}$ the complement of
$N(\q)$ in $I$, that is, $\overline{N(\q)}=I\setminus N(\q)$.

\begin{Lemma}\label{m2}
Let $\p$ and $\q$ be two distinct internal points. 

Assume that $q\equiv 1\pmod 4$.
\begin{itemize}
\item[(i)] If $\ell_{\p,\q}\in Pa_{\p}$ and 
$|\mathcal{U}_{\p,N(\q)}\cap C|$ is odd, then
$C=[\pi_k]$ for one $k$ or $C=D$.
\item[(ii)] If $\ell_{\p,\q}\in Se_{\p}$, then
$|\mathcal{U}_{\p, N(\q)}\cap C|$ is even.
\end{itemize}

Assume that $q\equiv 3\pmod 4$.
\begin{itemize}
\item[(iii)] If $\ell_{\p,\q}\in Pa_{\p}$, then 
$|\mathcal{U}_{\p,\overline{N(\q)}}\cap C|$ is even.

\item[(iv)] If $\ell_{\p,\q}\in Se_{\p}$ and
$|\mathcal{U}_{\p,\overline{N(\q)}}\cap C|$ is
odd, then $C=[\theta_i]$ for one $i$ or $C=D$.

\end{itemize}
\end{Lemma}

{\Proof} Without loss of generality, we can choose
$\p=(1,0,-\xi)$. Since $K=G_{\p}$ acts transitively
on both $Pa_{\p}$ and $Se_{\p}$, we may assume that 
$\q\not=\p$ is on either a special passant line 
$\ell_1=[1,0,\xi^{-1}]$ or a special secant line 
$\ell_2=[0,1,0]$ through $\q$.

{\bf Case I.} $\ell_1=\ell_{\p,\q}\in Pa_{\p}$.

In this case, $\q=(1,x,-\xi)$ for some $x\in\Ff_q$ 
with $u^2+\xi\in\Nsq_q$ and its internal neighbor is
$N(\q)=\{(1,u,-\xi)\mid u^2+\xi\in\Nsq_q\}\setminus\{(1,x,-\xi)\}$
by definition. As $\p\in N(\q)$, it is obvious that
$|\mathcal{U}_{\p,N(\q)}\cap D|$=1. Since the action 
of $K_{\q}$ on $Pa_{\q}$ has one orbit of length $1$, 
i.e. $\ell_1$, and all others are of length $2$, 
whose representatives form the set $\mathcal{R}_1$, 
we obtain that
\begin{equation}
\begin{array}{lllll}
|\mathcal{U}_{\p,N(\q)}\cap C| & = & 
\displaystyle\sum_{\ell\in Pa_\q}\displaystyle\sum_{\p_1\in I_{\ell}\setminus\{\q\}}
|\mathcal{U}_{\p,\p_1}\cap C|\\
{} & = & \displaystyle\sum_{\p_1\in I_{\ell_1}\setminus\{\q\}}|\mathcal{U}_{\p,\p_1}\cap C|+\
2\displaystyle\sum_{\ell\in \mathcal{R}}\displaystyle\sum_{\p_1\in I_{\ell}\setminus\{\q\}}
|\mathcal{U}_{\p,\p_1}\cap C|.\\
\end{array}
\end{equation}
Now let $\p_1=(1,u,-\xi)\in I_{\ell_1}\setminus\{\q\}$. 
Then the number of group elements that map $\p$ to $\p_1$ 
is determined by the quadruples $(a,b,c,d)$ which are the 
solutions to the following system of equations:
\begin{equation}\label{right2}
\begin{array}{cccc}
ab-cd\xi & = & u(a^2-c^2\xi)\\
b^2-d^2\xi & = & -\xi(a^2-c^2\xi)\\
a+d & = & s\\
ad-bc & = & 1.
\end{array}
\end{equation}
The first two equations in (\ref{right2}) yield 
$a^2-c^2\xi=A$ (or $-A$) where 
$A=\sqrt{\xi({u^2+\xi}^{-1})}$.

Now using 
$b^2-d^2\xi=\mp \xi A$, we obtain 
$$(b+c\xi)^2=s^2\xi-(2+A)\xi\;\;(\text{or}\;\; s^2\xi-(2-A)\xi).$$
If both $s^2\xi-(2+A)\xi$ and $s^2\xi-(2-A)\xi$ are 
squares, we set
$B_+=\sqrt{s^2\xi-(2+A)\xi}$ and $B_-=\sqrt{s^2\xi-(2-A)\xi}$,
then
$$a=\frac{1}{2s\xi}[s^2\xi-(B_\pm-2B_\pm \xi c)]
\;\;(\text{or}\;\;\frac{1}{2s\xi}[s^2\xi-(B_\pm+2B_\pm \xi c)])$$
and 
$$d=\frac{1}{2s\xi}[s^2\xi+(B_\pm-2B_\pm \xi c)]
\;\;(\text{or}\;\;\frac{1}{2s\xi}[s^2\xi+(B_\pm+2B_\pm \xi c)]);$$
combining with the last two equations of (\ref{right2}),
we have
\begin{equation}\label{dis} 
(\xi-\frac{B_\pm^2}{s^2})c^2+(\frac{B_\pm^3}{s^2\xi}-B_\pm)c+
(\frac{s^2}{4}-\frac{B_\pm^4}{4s^2\xi^2}-1)=0
\end{equation}
or
\begin{equation}\label{dis_1} 
(\xi-\frac{B_\pm^2}{s^2})c^2-(\frac{B_\pm^3}{s^2\xi}-B_\pm)c+
(\frac{s^2}{4}-\frac{B_\pm^4}{4s^2\xi^2}-1)=0.
\end{equation}
The discriminant of (\ref{dis}) or (\ref{dis_1}) is
$$\Delta=(1-\frac{B_\pm^2}{s^2\xi})(B_\pm^2-s^2\xi+4\xi)=
\frac{4\xi u^2}{s^2(u^2+\xi)}\in\Sq_q.$$
Consequently, the equations in (\ref{right2}) have $8$ solutions
and yield $4$ different group elelments.

If one of $s^2\xi-(2+A)\xi$ and $s^2\xi-(2-A)\xi$ is a square
and the other is nonsquare, similar arguments as above give that
the equations in (\ref{right2}) have $4$ solutions and produce 
$2$ different group element.  

If one of $s^2\xi-(2+A)\xi$ and $s^2\xi-(2-A)\xi$ is zero, then
$s^2$ is one of $2+A$ and $2-A$; and moreover it is one of $\pi_k$ 
for $1\le k\le \frac{q-1}{4}$ since 
$(2+A)(2-A)=\frac{4u^2}{u^2+\xi}\in\Nsq_q$ and $-1\in\Sq_q$. 
Consequently, the equations in (\ref{right2}) yield either $1$ or 
$3$ group elements in $[s^2]$.

Therefore, if $|\mathcal{U}_{\p,N(\q)}\cap C|$ is odd,
then $C=D$ or $[\pi_k]$ for one $k$.
\vskip.1in
{\bf Case II.} $\ell_2=\ell_{\p,\q}\in Se_{\p}$ and $\q\notin\p^\perp$.

Then $\q=(1,0,-y)$ for $y\notin\Nsq_q$ and $y\not=\pm\xi$. 
From the proof of Case II in Lemma~\ref{m1}, we have that
$K_{\q}=\{{\bf d}(1,1,1), {\bf ad}(-1,1,-1)\}$, and among the
orbits of $K_\q$ on $Pa_\p$, $K_\q$ has
only one orbit of length $1$, that is, $\ell_4=[1,0,y^{-1}]$;
and all other orbits are of length $2$ whose representatives
form the set $\mathcal{R}$. Since 
$|\mathcal{U}_{\p,I_{\ell_i}}\cap C|
=|\mathcal{U}_{\p,I_{\ell_j}}\cap C|$ 
where $\ell_i$, $\ell_j\in Pa_\p$ and $\ell_j=\ell_i^g$ for $g\in K_\q$, 
we obtain that
\begin{equation}\label{summ1}
\begin{array}{lllll}
|\mathcal{U}_{\p,N(\q)}\cap C| & = & \displaystyle\sum_{\ell\in Pa_{\q}}
\displaystyle\sum_{\p_1\in I_{\ell}\setminus\{\q\}}|\mathcal{U}_{\p,\p_1}\cap C|\\
{} & = & \displaystyle\sum_{\p_1\in I_{\ell_4}\setminus\{\q\}}|\mathcal{U}_{\p,\p_1}\cap C|+2
\displaystyle\sum_{\ell\in \mathcal{R}}\displaystyle\sum_{\p_1\in I_{\ell}\setminus\{\q\}}
|\mathcal{U}_{\p,\p_1}\cap C|.\\
\end{array}
\end{equation}
Moreover, since the orbits of $K_\q$ on $I_{\ell_4}\setminus\{\q\}$,
whose representatives form the set $\mathcal{R}_1$, are of length 
$2$ and $|\mathcal{U}_{\p,\p_1}\cap C|=|\mathcal{U}_{\p,\p_2}\cap C|$ 
for $\p_2=\p_1^g$, the first term of the last expression in 
(\ref{summ1}) can be rewritten as 
\begin{displaymath}
2\displaystyle\sum_{\p_1\in \mathcal{R}_1}|\mathcal{U}_{\p,\p_1}\cap C|.
\end{displaymath}
So $|\mathcal{U}_{\p,N(\q)}\cap C|$ is even in this case.
\vskip.1in
{\bf Case III.} $\p=\ell_2\cap \p^\perp$.

In this case, we have $\q=(1,0,\xi)$. Among the orbits of $K_\q$ on
$Pa_\p$, only one has length $1$, i.e. $\p^\perp$. Moreover, all the
orbits of $K_\q$ on $I_{\p^\perp}\setminus\{\q\}$ are of length $2$. 
Hence $|\mathcal{U}_{\p,N(\q)}\cap C|$ is even.

The case when $q\equiv 3\pmod 4$ can be established in the same
way and we eliminate the detail. \QED



\section{Group Algebra $FH$}
\subsection{2-Blocks of H}
In this section we recall several results on the $2$-blocks of $H\cong \PSL(2,q)$.
We refer the reader to {\rm \cite{gabriel}} or {\rm
{\cite{brauer}}} for a general introduction on this subject.

Let $\mathbf{R}$ be the ring of algebraic integers in the complex
field $\Cc$. We choose a maximal ideal $\mathbf{M}$ of $\mathbf{R}$
containing $2\mathbf{R}$. Let $F=\mathbf{R}/\mathbf{M}$ be the
residue field of characteristic $2$, and let $* :
\mathbf{R}\rightarrow F$ be the natural ring homomorphism. Define
\begin{equation}
\begin{array}{llll}\label{ring_s}
\mathbf{S}& =& \{\frac{r}{s}\mid r\in\mathbf{R},\;s\in\mathbf{R}\setminus\mathbf{M}\}.\\
\end{array}
\end{equation}
Then it is clear that the map $* : \mathbf{S}\rightarrow F$ defined
by $(\frac{r}{s})^* = r^*(s^*)^{-1}$ is a ring homomorphism with
kernel $\mathcal{P} = \{\frac{r}{s}\mid
r\in\mathbf{M},\;s\in\mathbf{R}\setminus\mathbf{M}\}$. In the rest
of this article, $F$ will always be the field of characteristic $2$
constructed as above. Note that $F$ is an algebraic closure of
$\Ff_2$.

As a convention, we use $\textrm{Irr}(H)$ and $\textrm{IBr}(H)$
to denote the set of irreducible ordinary characters and the set 
of irreducible Brauer characters of $H$, respectively. 

In the following, the irreducible characters of $H$ 
are classified according to the character tables of 
$\PSL(2,q)$ displayed in Appendix.

\begin{Lemma}$($\cite{jordan}, \cite{GK}, \cite{schur}$)$
The irreducible ordinary characters of $H$ are:
\begin{itemize}
\item[(i)] $1=\chi_0$, $\gamma$, $\chi_1$, ..., $\chi_{\frac{q-1}{4}}$,
$\beta_1$, $\beta_2$, $\phi_1$, ..., $\phi_{\frac{q-5}{4}}$ if 
$q\equiv 1\pmod 4$, where $1=\chi_0$ is the trivial character,
$\gamma$ is the character of degree $q$, $\chi_s$ for $1\le s\le \frac{q-1}{4}$ 
are the characters of degree $q-1$, $\phi_r$ for $1\le r\le \frac{q-5}{4}$
are the characters of degree $q+1$, and $\beta_i$ for $i=1$, $2$ are the
characters of degree $\frac{q+1}{2}$;
\item[(ii)]  $1=\chi_0$, $\chi_1$, ..., $\chi_{\frac{q-3}{4}}$,
$\beta_1$, $\eta_2$, $\eta_1$, ..., $\phi_{\frac{q-3}{4}}$ if 
$q\equiv 3\pmod 4$, where $1=\chi_0$ is the trivial character,
$\gamma$ is the character of degree $q$, $\chi_s$ for $1\le s\le \frac{q-1}{3}$
are the characters of degree $q-1$, $\phi_r$ for $1\le r\le \frac{q-3}{4}$
are the characters of degree $q+1$, and $\eta_i$ for $i=1$, $2$ are the
characters of degree $\frac{q-1}{2}$;
\end{itemize}
\end{Lemma}
The following lemma illumstrates how the irreducible ordinary 
characters of $H$ are partitioned into $2$-blocks.
\begin{Lemma}\cite[Lemma 4.1]{swx}\label{blocks}
First assume that $q\equiv 1 \pmod 4$ and $q-1 = m2^n$, where 
$2\nmid m$.
\begin{itemize}

\item[(i)] The principal block $B_0$ of $H$ contains $2^{n-2}+3$ 
irreducible characters
$\chi_0 = 1$, $\gamma$, $\beta_1$, $\beta_2$, $\phi_{i_1}$, ..., 
$\phi_{i_{(2^{n-2}-1)}}$, 
where $\chi_0=1$ is the trivial character of $H$, $\gamma$ is 
the irreducible character of degree $q$ of $H$,  $\beta_1$ and 
$\beta_2$ are the two irreducible characters of degree 
$\frac{q+1}{2}$, and $\phi_{i_k}$ for $1\le k\le 2^{n-2}-1$ 
are distinct irreducible characters of degree $q+1$ of $H$.

\item[(ii)]  $H$ has $\frac{q-1}{4}$ blocks $B_s$ of defect 
$0$ for $1\le s \le \frac{q-1}{4}$, each of which contains 
an irreducible ordinary character $\chi_s$ of degree $q-1$.

\item[(iii)] If $m\ge 3$, then $H$ has $\frac{m-1}{2}$ blocks 
$B_t^{'}$ of defect $n-1$ for $1\le t\le \frac{m-1}{2}$, 
each of which contains $2^{n-1}$ irreducible ordinary characters 
$\phi_{t_i}$ for $1\le i \le 2^{n-1}$.

\end{itemize}

Now assume that $q\equiv 3\pmod 4$ and $q+1=m2^n$, where 
$2\nmid m$ .
\begin{itemize}

\item[(iv)] The principal block $B_0$ of $H$ contains 
$2^{n-2}+3$ irreducible characters
$\chi_0 = 1$, $\gamma$, $\eta_1$, $\eta_2$, $\chi_{i_1}$, 
..., $\chi_{i_{(2^{n-2}-1)}}$,
where $\chi_0=1$ is the trivial character of $H$, $\gamma$ 
is the irreducible character of degree $q$ of $H$, $\eta_1$ 
and $\eta_2$ are the two irreducible characters of degree
$\frac{q-1}{2}$, and $\chi_{i_k}$ for $1\le k\le 2^{n-2}-1$ 
are distinct irreducible characters of degree $q-1$ of $H$.

\item[(v)] $H$ has $\frac{q-3}{4}$ blocks $B_r$ of defect 
$0$ for $1\le r\le \frac{q-3}{4}$, each of which contains 
an irreducible ordinary character $\phi_r$ of degree $q+1$.

\item[(vi)] If $m\ge 3$, then $H$ has $\frac{m-1}{2}$ blocks 
$B_t^{'}$ of defect $n-1$ for $1\le t\le \frac{m-1}{2}$, each 
of which contains $2^{n-1}$ irreducible ordinary characters 
$\chi_{t_i}$ for $1\le i \le 2^{n-1}$.

\end{itemize}

Moreover, the above blocks form all the $2$-blocks of $H$.
\end{Lemma}
\begin{Remark}
Parts $($i$)$ and $($iv$)$ are from Theorem 1.3 in 
{\rm\cite{landrock}} and their proofs can be found in 
Chapter 7 of {III} in {\rm \cite{brauer}}. Parts $($ii$)$ 
and $($v$)$ are special cases of Theorem 3.18 in 
{\rm \cite{gabriel}}. Parts $($iii$)$ and $($vi$)$ are
proved in Sections {II} and {VIII} of {\rm \cite{burkhardt}}.
\end{Remark}


\subsection{Block Idempotents}

Let $Bl(H)$ be the set of $2$-blocks of $H$. If 
$B\in Bl(H)$, we write
$$f_B = \displaystyle\sum_{\chi\in Irr(B)}e_{\chi},$$
where $e_{\chi}=\frac{\chi(1)}{|H|}\sum_{g\in H} \chi(g^{-1})g$ 
is a central primitive idempotent of $\mathbf{Z}(\Cc H)$ 
and $Irr(B)=Irr(H)\cap B$. For future use, we define 
$IBr(B)=IBr(H)\cap B$. Since $f_B$ is an element of 
$\mathbf{Z}(\Cc H)$, we may write
\begin{displaymath}
\begin{array}{lllll}
f_B& =& \displaystyle\sum_{C\in cl(H)}f_B(\widehat{C})\widehat{C},\\
\end{array}
\end{displaymath}
where $cl(H)$ is the set of conjugacy classes of $H$, 
$\widehat{C}$ is the sum of elements in the class $C$, 
and
\begin{equation}\label{id}
\begin{array}{llll}
f_B(\widehat{C})& = &\frac{1}{|H|}\displaystyle\sum_{\chi\in Irr(B)}\chi(1)\chi(x_C^{-1})
\end{array}
\end{equation}
with a fixed element $x_C\in C$.
\begin{Theorem}\label{osima}
Let $B\in Bl(H)$. Then $f_B\in \mathbf{Z}(\mathbf{S}H)$. 
In other words, $f_B(\widehat{C})\in \mathbf{S}$ for each 
block of $H$.

\end{Theorem}
{\Proof} It follows from Corollary 3.8 in {\rm \cite{gabriel}}. \QED

We extend the ring homomorphism $*: \mathbf{S}\rightarrow F$ 
to a ring homomorphism $*:\mathbf{S}H\rightarrow FH$ by setting
$(\sum_{g\in H} s_g g)^*= \sum_{g\in H} s_g^* g$. Note that $*$ 
maps $\mathbf{Z}(\mathbf{S}H)$ onto $\mathbf{Z}(FH)$ via 
$(\sum_{C\in cl(H)}s_C \widehat{C})^*$ = 
$\sum_{C\in cl(H)} s_C^* \widehat{C}$.
Now we define
$$e_B = (f_B)^* \in \mathbf{Z}(FH),$$
which is the {\it block idempotent} of $B$. Note that $e_B
e_{B^{'}}= \delta_{B B^{'}}e_B$ for $B$, $B^{'}\in Bl(H)$, where
$\delta_{B B^{'}}$ equals 1 if $B=B'$, 0 otherwise. Also
$1=\sum_{B\in Bl(H)}e_B$.

All the block idempotents of the $2$-blocks of $H$ are given in the
following lemma; see \cite{swx} for the detailed calculations.

\begin{Lemma}\cite[Lemma 4.4]{swx}\label{expression}
First assume that $q\equiv 1 \pmod 4$ and $q-1= m2^n$ with $2\nmid m$.
\begin{itemize}

\item[1.] Let $B_0$ be the principal block of $H$. Then
$($a$)$ $e_{B_0}(\widehat{D})=1$;
$($b$)$ $e_{B_0}(\widehat{F^+})= e_{B_0}(\widehat{F^-})\in F$;
$($c$)$ $e_{B_0}(\widehat{[\theta_i]})\in F$, $e_{B_0}(\widehat{[0]})=0$;
$($d$)$ $e_{B_0}(\widehat{[\pi_k]})=1$.
\item[2.] Let $B_s$ be any block of defect $0$ of $H$. Then
$($a$)$ $e_{B_s}(\widehat{D})=0$;
$($b$)$ $e_{B_s}(\widehat{F^+})= e_{B_s}(\widehat{F^-})=1$;
$($c$)$ $e_{B_s}(\widehat{[0]})= e_{B_s}(\widehat{[\theta_i]})= 0$;
$($d$)$ $e_{B_s}(\widehat{[\pi_k]})\in F$.
\item[3.] Suppose $m\ge 3$ and let $B_t^{'}$ be any block of defect $n-1$ of $H$. Then
$($a$)$ $e_{B_t^{'}}(\widehat{D})=0$;
$($b$)$ $e_{B_t^{'}}(\widehat{F^+})= e_{B_t^{'}}(\widehat{F^-})=1$;
$($c$)$ $e_{B_t^{'}}(\widehat{[\theta_i]})\in F$, $e_{B_t^{'}}(\widehat{[0]})=0$;
$($d$)$ $e_{B_t^{'}}(\widehat{[\pi_k]}) = 0$.
\end{itemize}
Now assume that $q\equiv 3 \pmod 4$. Suppose that $q+1=m2^n$ with $2\nmid m$.
\begin{itemize}

\item[4.] Let $B_0$ be the principal block of $H$. Then
$($a$)$ $e_{B_0}(\widehat{D})=1$;
$($b$)$ $e_{B_0}(\widehat{F^+})= e_{B_0}(\widehat{F^-})\in F$;
$($c$)$ $e_{B_0}(\widehat{[\theta_i]})=1$;
$($d$)$ $e_{B_0}(\widehat{[0]})=0$, $e_{B_0}(\widehat{[\pi_k]})\in F$.
\item[5.] Let $B_r$ be any block of defect $0$ of $H$. Then
$($a$)$ $e_{B_r}(\widehat{D})=0$;
$($b$)$ $e_{B_r}(\widehat{F^+})= e_{B_r}(\widehat{F^-})=1$;
$($c$)$ $e_{B_r}(\widehat{[0]})= e_{B_r}(\widehat{[\pi_k]}) = 0$;
$($d$)$ $e_{B_r}(\widehat{[\theta_i]})\in F$.
\item[6.] Suppose that $m\ge 3$ and let $B_t^{'}$ be any block of defect $n-1$ of $H$. Then
$($a$)$ $e_{B_t^{'}}(\widehat{D})=0$;
$($b$)$ $e_{B_t^{'}}(\widehat{F^+})= e_{B_t^{'}}(\widehat{F^-})=1$;
$($c$)$ $e_{B_t^{'}}(\widehat{[\theta_i]})= 0$;
$($d$)$ $e_{B_t^{'}}(\widehat{[0]})=0$, $e_{B_t^{'}}(\widehat{[\pi_k]})\in F$.
\end{itemize}
\end{Lemma}

Let $M$ be an $\mathbf{S}H$-module. We denote the 
reduction $M/\mathcal{P}M$, which is an $FH$-module, 
by $\overline{M}$. Then the following lemma is apparent.
\begin{Lemma}\label{reduction}
Let $M$ be an $\mathbf{S}H$-module and $B\in Bl(H)$. 
Using the above notation, we have \
$$\overline{M\cdot f_b} = \overline{M}\cdot e_B,$$
i.e. reduction commutes with projection onto a block 
$B$.
\end{Lemma}

\section{Linear Maps and Their Matrices}

Let $F$ be the algebraic closure of $\Ff_2$ defined 
in Section 4. Recall that for $\p\in I$, $N(\p)$ is 
the set of external points on the passant lines 
through $\p$ with $\p$ included or excluded 
accordingly as $q\equiv 3\pmod 4$ or $q\equiv 1\pmod 4$.
We define $\mathbf{D}$ to be the incidence matrix 
of $N(\p)$ ($\p\in I$) and $I$. Correspondingly, 
the rows of $\mathbf{D}$ can be viewed as the 
characteristic vectors of $N(\p)$. In the following, 
we always regard both $\D$ and $\A$ as matrices over $F$.
Moreover, it is clear that $\D=\A^2+\mathbf{I}$,
where $\mathbf{I}$ is the identity matrix of proper
size.

\begin{Definition}
For $W\subseteq\in I$, we define $\mathcal{C}_{W}$ 
to be the row characteristic vector of $W$ with 
respect to $I$, i.e. $\mathcal{C}_{W}$ is a $0$-$1$ 
row vector of length $|I|$ with entries indexed 
by the internal points and the entry of $\mathcal{C}_{W}$ 
is $1$ if and only if the point indexing the entry
is in $W$. If $W=\{\p\}$, as a convention, we write
$\mathcal{C}_W$ as $\mathcal{C}_\p$.
\end{Definition}

Let $k$ be the complex field $\Cc$, the algebraic 
closure $F$ of $\Ff_2$, or the ring $\mathbf{S}$ 
in~(\ref{ring_s}). Let $k^I$ be the free $k$-module 
with the base 
$\{\mathcal{C}_{\p}\mid \p\in I\}$, respectively. 
If we extend the action of $H$ on the basis elements
of $k^I$, which is defined by
$\mathcal{C}_{\q}\cdot h=\mathcal{C}_{\q^h}$ for 
$\p \in I$ and $h\in H$, linearly to $k^I$, then 
$k^I$ is a $kH$-permutation modules. Since $H$ is 
transitive on $I$, we have 
$$k^I=\textrm{Ind}_K^H(1_k),$$
where $K$ is the stabilizer of an internal point 
in $H$ and $\textrm{Ind}_K^H(1_k)$ is the $kH$-module 
induced from $1_k$.

The decomposition of $1\uparrow_K^H$, the character of 
$\textrm{Ind}_K^H(1_k)$, into a sum of the irreducible 
ordinary characters of $H$ is given as follows.

\begin{Lemma}\label{decomposition_1}
Let $K$ be the stabilizer of an internal
point in $H$. 

Assume that $q\equiv 1\pmod 4$. Let $\chi_s$,
$1\le s \le \frac{q-1}{4}$, be the irreducible
characters of degree $q-1$, $\phi_r$, 
$1\le r \le \frac{q-5}{4}$,
irreducible characters of degree $q+1$, $\gamma$
the irreducible character of degree $q$, and $\beta_j$, 
$1\le j\le 2$, irreducible characters of degree 
$\frac{q+1}{2}$.

\begin{itemize}
\item[(i)] If $q\equiv 1 \pmod 8$, then $$1_K\uparrow_{K}^H 
=1_H+\displaystyle\sum_{s=1}^{(q-1)/4}\chi_s + \gamma +
\beta_1+\beta_2 + \displaystyle\sum_{j=1}^{(q-9)/4}\phi_{r_j},$$
where $\phi_{r_j}$, $1\le j\le\frac{q-9}{4}$, may not be 
distinct.
\item[(ii)] If $q\equiv 5\pmod 8$, then $$1_K\uparrow_{K}^H =
1_H+\displaystyle\sum_{s=1}^{(q-1)/4}\chi_s + \gamma +
\displaystyle\sum_{j=1}^{(q-5)/4}\phi_{r_j},$$ where
$\phi_{r_j}$, $1\le j\le\frac{q-5}{4}$, may not be 
distinct.
\end{itemize}

Next assume that $q\equiv 3\pmod 4$. Let $\chi_s$, 
$1\le s \le\frac{q-3}{4}$, be the irreducible
characters of degree $q-1$, $\phi_r$, 
$1\le r \le \frac{q-3}{4}$, the irreducible
characters of degree $q+1$, $\gamma$ the irreducible
character of degree $q$, and $\eta_j$, $1\le j\le 2$, 
the irreducible characters of degree $\frac{q-1}{2}$.

\begin{itemize}
\item[(iii)] If $q\equiv 3\pmod 8$, then 
$$1_K\uparrow_{K}^H= 1_H+\displaystyle\sum_{r=1}^{(q-3)/4}
\phi_r+\eta_1+\eta_2+\displaystyle\sum_{j=1}^{(q-3)/4}
\chi_{s_j},$$
where $\chi_{s_j}$, $1\le j\le\frac{q-3}{4}$, may not 
be distinct.

\item[(iv)] If $q\equiv 7\pmod 8$, then 
$$1_K\uparrow_{K}^H=1_H+\displaystyle\sum_{r=1}^{(q-3)/4}
\phi_r+\displaystyle\sum_{j=1}^{(q+1)/4}\chi_{s_j},$$ 
where $\chi_{s_j}$, $1\le j\le\frac{q+1}{4}$, may not 
be distinct.
\end{itemize}
\end{Lemma}
{\Proof} We provide the proof for the case when 
$q\equiv 1\pmod 4$.

Let $1_H$ be the trivial character of $H$. By the 
Frobenius reciprocity \cite{frob},
$$\left\langle 1_K\uparrow_{K}^H, 1_H\right\rangle_H=
\left\langle 1_K, 1_H\downarrow_K^H\right\rangle_K = 1.$$

Let $\chi_s$ be an irreducible character of degree $q-1$
of $H$, where $1\le s \le \frac{q-1}{4}$. We denote the
number of elements of $K$ lying in the class $[\pi_k]$ by
$d_k$. Then $d_k=2$
by Lemma~\ref{Kintersection}(iii), and so
\begin{displaymath}
\begin{array}{lllll}
\left\langle1_K\uparrow_K^H,\chi_s\right\rangle_H & = &
\left\langle1_K, \chi_s\downarrow_K^H\right\rangle_K
& = &\frac{1}{|K|}\displaystyle\sum_{g\in K}\chi_s
\downarrow_{K}^H(g)\\
{} & {} & {} &=&\frac{1}{q+1}[(1)(q-1)+2\displaystyle\sum_{k=1}^{(q-1)/4}
(-\delta^{(2k)s}-\delta^{-(2k)s})]\\
{} & {} &{} &=& 1,
\end{array}
\end{displaymath}
where 
\begin{displaymath}
\begin{array}{lll}
\displaystyle\sum_{k=1}^{(q-1)/4}(-\delta^{(2k)s}-\delta^{-(2k)s})&=&
-(1+\delta^{2s}+(\delta^{2s})^2+\cdots+(\delta^{2s})^{(q-1)/2}-1)\\
{} & = & -\frac{1-\delta^{(q+1)s}}{1-\delta^{2s}}+1\\
 {} & = & 1
\end{array}
\end{displaymath}
since $\delta^{q+1}=1$.

Let $\gamma$ be the irreducible character of degree $q$ of $H$.
Then
\begin{displaymath}
\begin{array}{llllll}
\left\langle 1_K\uparrow_K^H,\gamma \right\rangle_H &=& 
\left\langle 1_K, \gamma\downarrow_K^H\right\rangle_K 
&=& \frac{1}{|K|}\displaystyle\sum_{g\in K} \gamma
\downarrow_K^H(g)\\
{}&{}&{}&=& \frac{1}{q+1}[(1)(q)+(2)(-1)(\frac{q-1}{4})+
(1)(\frac{q+1}{2})]\\
{} & {} & {} &=& 1.
\end{array}
\end{displaymath}

Let $\beta_j$ be any irreducible character of degree
$\frac{q+1}{2}$ of $H$. Then
\begin{equation}\label{charsum1}
\begin{array}{llll}
\left\langle1_K\uparrow_K^H, \beta_j\right\rangle_H &=&
\frac{1}{|K|}\displaystyle\sum_{g\in K}\beta_j\downarrow_K^H(g)\\
{} & 
= &\frac{1}{q+1}[(1)(\frac{q+1}{2})+(2)
(\frac{q-1}{4})(0)+(\frac{q+1}{2})(-1)^{(q-1)/4}].
\end{array}
\end{equation}
Consequently, if $q\equiv 1\pmod 8$, then
$(-1)^{\frac{q-1}{4}}=1$, and so 
$\left\langle1_K\uparrow_K^H, \beta_j\right\rangle_H=1$;
otherwise, $(-1)^{\frac{q-1}{4}}=-1$, and so 
$\left\langle1_K\uparrow_K^H, \beta_j\right\rangle_H=0$.

Since the sum of the degrees of $1$, $\chi_s$, 
$\gamma$, and $\beta_j$ is less than the degree 
of $1\uparrow_K^H$ and only the irreducible 
characters of degree $q+1$ of $H$ have not been
taken into account yet, we see that all the 
irreducible constituents of
$$1_K\uparrow_K^H - 1_H - \displaystyle\sum_{s=1}^{(q-1)/4}
\chi_s-\gamma-\beta_1-\beta_2\;\;
\textrm{or}
\;\;1_K\uparrow_K^H-1_H-\displaystyle\sum_{s=1}^{(q-1)/4}
\chi_s-\gamma$$
must have degree $q+1$. \QED

\begin{Corollary}\label{char11}
Using the above notation,  
\begin{itemize}
\item[(i)] if $q\equiv 1\pmod 4$, then the character of 
$\mathrm{Ind}_K^H (1_{\Cc})\cdot f_{B_s}$ is $\chi_s$ 
for each block $B_s$ of defect $0$;
\item[(ii)] if $q\equiv 3\pmod 4$, then the character of 
$\mathrm{Ind}_K^H (1_{\Cc})\cdot f_{B_r}$ is $\phi_r$ 
for each block $B_r$ of defect $0$. 
\end{itemize}
\end{Corollary}
{\Proof} The corollary follows from Lemma~\ref{blocks} 
and Lemma~\ref{decomposition_1}.
\QED

Since $H$ preserves incidence, it is obvious that,
for $\p\in I$ and $h\in H$,
$$ 
\mathcal{C}_{N(\p)}\cdot h = \mathcal{C}_{N(\p^h)}.
$$

In the rest of the article, we always view 
$\mathcal{C}_{\p}$ as a vector over $F$. Consider the 
maps $\phi$ and $\mu$ from $F^I$ to $F^I$ 
defined by extending
$$\mathcal{C}_{\p}\mapsto \mathcal{C}_{\p^\perp}, 
\mathcal{C}_{\p}\mapsto \mathcal{C}_{N(\p)}$$
linearly to $F^I$, respectively. Then it is clear that 
as $F$-linear maps, the marices of $\phi$ and $\mu$, 
are $\A$ and $\D$, respectively,
and for ${\bf x}\in F^I$, $\phi({\bf x})={\bf x}\A$
and $\mu({\bf x})={\bf x}\D$. Moreover, we have the 
following result since $H$ is transitive on $I$ and
preserves incidence:
\begin{Lemma}\label{hom}
The maps $\phi$ and $\mu$ are both $FH$-module 
homomorphisms from $F^I$ to $F^I$.
\end{Lemma}

We will always use ${\bf 0}$ and $\hat{{\bf 0}}$ 
to denote the all-zero row vector of length $|I|$ 
and the all-zero matrix of size $|I|\times |I|$, 
respectively; and we denote by $\hat{\mathbf{J}}$
and $\mathbf{J}$ the all-one row vector of length $|I|$
and the all-one matrix of size $|I|\times |I|$.
\begin{Proposition}\label{sum}
As $FH$-modules, $F^I=\Ima(\phi)\oplus \Ker(\phi)$, 
where $\Ima(\phi)$ and $\Ker(\phi)$ are the image 
and kernel of $\phi$, respectively.
\end{Proposition}
{\Proof} It is clear that 
$\Ker(\phi)\subseteq \Ker(\phi^2)$. If
${\bf x}\in \Ker(\phi^2)$, then ${\bf x}\in \Ker(\phi)$ 
since
$$\phi({\bf x})=\phi^3({\bf x})=\phi(\phi^2({\bf x}))={\bf 0}.$$
Therefore, $\Ker(\phi^2)=\Ker(\phi)$. Furthermore, 
since 
$\Ker(\phi)\subseteq \Ker(\phi^2)\subseteq \Ker(\phi^3)\subseteq\cdots$,
we have $\Ker(\phi^i)=\Ker(\phi)$ for $i\ge 2$.
Applying the Fitting decomposition theorem \cite[p. 285]{TYY} 
to the operator $\phi$, we can find an $i$ such that 
$F^E=\Ker(\phi^i)\oplus \Ima(\phi^i)$. From the above 
discussions, we must have 
$F^E=\Ker(\phi)\oplus \Ima(\phi)$.
\QED


\begin{Corollary}\label{directsum}
As $FH$-modules, $\Ind_K^H(1_F)\cong \Ker(\phi)\oplus \Ima(\phi)$.
\end{Corollary}
{\Proof} The conclusion follows immediately from
Proposition~\ref{sum} and the fact that $\Ind_K^H(1_F)\cong
F^E$. \QED

Using the above notation, we set $\mathbf{C}=\D+\jj$, where $\jj$
is the all-one matrix of proper size.
Then the matrix $\mathbf{C}$ can be viewed as the 
incidence matrix of $\overline{N(\p)}$ ($\p\in I$) and $I$, 
and so $\mathcal{C}_\p\C=\mathcal{C}_{\overline{N(\p)}}$. 

Let $\mu_2$ be the $FH$-homomorphism 
from $F^I$ to $F^I$ whose matrix with respect to the natural 
basis is $\mathbf{C}$. 
The following proposition is clear.
\begin{Proposition}\label{proker}
Using the above notation, we have
$\Ker(\phi)=\Ima(\mu)$.
\end{Proposition}
Furthermore, we have the following decomposition of $\Ker(\phi)$.
\begin{Lemma}\label{separate}
Assume that $q\equiv 3\pmod 4$. Then as $FH$-modules, we have
$\Ker(\phi)=\langle\hat{\mathbf{J}}\rangle\oplus \Ima(\mu_2)$, 
where
$\langle\hat{\mathbf{J}}\rangle$ is the trivial $FH$-module
generated by $\hat{\bf J}$.
\end{Lemma}
{\Proof} Let ${\bf y}\in \langle\hat{\mathbf{J}}\rangle\cap
\Ima(\mu_2)$. Then ${\bf y}=\mu_2({\bf x})=\lambda\hat{\mathbf{J}}$
for some $\lambda\in F$ and ${\bf x}\in F^I$. Or equivalently,
we have $\mu_2({\bf x})={\bf x}\C={\bf x}(\A^2+{\bf I}+{\bf J})
=\lambda\hat{\mathbf{J}}$.
Note that
${\jj}^2={\jj}$ and $\hat{\bf {J}}{\jj}=\hat{\bf J}$ since 
$2\nmid |I|$ when $q\equiv 3 \pmod 4$. Moreover, 
$\A^2\jj=\hat{\bf 0}$ as each row of $\A^2$ , viewed as the
characteristic vector of $\widehat{N(\p)}$, has an even number of 
$1$s. Consequently, 
$$\lambda\hat{\bf J}=\lambda\hat{\bf J}\jj={\bf x}
(\mathbf{A}^2+{\ii}+{\jj}){\jj} = {\bf x}(\mathbf{A}^2{\jj}+{\ii}{\jj}+{\jj}^2)
={\bf x}(\hat{{\bf 0}}+ {\jj} + {\jj}) ={\bf 0}.$$
It follows that $\lambda=0$. Therefore, we must have
$\langle\hat{\mathbf{J}}\rangle\cap \Ima(\mu_2)={\bf 0}$.

It is obvious that $\langle\hat{\mathbf{J}}\rangle +
\Ima(\mu_2)\subseteq \Ker(\phi)$. Let ${\bf x}\in \Ker(\phi)$. 
Then ${\bf x}={\bf y}(\mathbf{A}^2+{\ii})$ for some 
${\bf y}\in F^I$. Since 
${\bf y}{\jj}=\langle{\bf y}, \hat{{\bf J}}\rangle\hat{{\bf J}}$,
we obtain that 
${\bf x}={\bf y}(\mathbf{A}^2+{\ii}+ {\jj})+\langle{\bf y}, 
\hat{{\bf J}}\rangle\hat{{\bf J}}$,
where $\langle{\bf y}, \hat{{\bf J}}\rangle$ is the standard 
inner product of the vectors ${\bf y}$ and $\hat{{\bf J}}$.
Hence ${\bf x}\in\langle\hat{\mathbf{J}}\rangle + \Ima(\mu_2)$ 
and so 
$\Ker(\phi)=\langle\hat{\mathbf{J}}\rangle\oplus \Ima(\mu_2)$. \QED



\section{Statement and Proof of Main Theorem}
The main theorem is stated as follows.

\begin{Theorem}\label{main}
Let $\Ker(\phi)$ be defined as above. As $FH$-modules,
\begin{itemize}
\item[(i)] if $q\equiv 1\pmod 4$, then

$$\Ker(\phi)=\displaystyle\bigoplus_{s=1}^{(q-1)/4}M_s,$$
where $M_s$ for $1\le s\le \frac{q-1}{4}$ are pairwise 
non-isomorphic simple $FH$-modules of dimension $q-1$;

\item[(ii)] if $q\equiv 3\pmod 4$, then

$$\Ker(\phi)=\langle\hat{\mathbf{J}}\rangle
\oplus (\displaystyle\bigoplus_{r=1}^{(q-3)/4}M_r),$$
where $M_r$ for $1\le s\le \frac{q-3}{4}$ are pairwise 
non-isomorphic simple $FH$-modules of dimension $q+1$ 
and $\langle\hat{\mathbf{J}}\rangle$ is the trivial
$FH$-module generated by the all-one column vector of 
length $|I|$.

\end{itemize}

\end{Theorem}
To prove the main theorem, we need refer to the following 
lemma.

\begin{Lemma}\label{y4}
Let  $q-1=2^n m$ or $q+1=2^n m$
with $2\nmid m$ accordingly as $q\equiv 1\pmod 4$ 
or $q\equiv 3\pmod 4$. Using the above notation, 
\begin{itemize}
\item[(i)] if $q\equiv 1\pmod 4$, then $\Ker(\phi)\cdot e_{B_0}
={\bf 0}$, $\Ima(\phi)\cdot e_{B_s}={\bf 0}$ for $1\le s\le 
\frac{q-1}{4}$, and $\Ker(\phi)\cdot e_{B_t^{'}}={\bf 0}$ for 
$m\ge 3$ and $1\le t \le \frac{m-1}{2}$;

\item[(ii)] if $q\equiv 3\pmod 4$, then 
$\Ima(\mu_2)\cdot e_{B_0}={\bf 0}$, 
$\Ima(\phi)\cdot e_{B_r}={\bf 0}$ for $1 \le r \le \frac{q-3}{4}$, and 
$\Ima(\mu_2)\cdot e_{B_t^{'}}={\bf 0}$ for $m\ge 3$ and 
$1\le t\le \frac{m-1}{2}$.
\end{itemize}

\end{Lemma}

{\Proof} It is clear that $\Ima(\phi)$, $\Ker(\phi)$, and 
$\Ima(\mu_2)$ are generated by 
$$\{\mathcal{C}_{\p^\perp}\mid \p\in I\},\;\; 
\{\mathcal{C}_{N(\p)}\mid \p\in I\},\;\;\text{and}\;\;
\{\mathcal{C}_{\overline{N(\p)}}\mid \p\in I\}$$
over $F$, respectively. Now let $B\in Bl(H)$. Since
\begin{equation*}
\begin{array}{llll}
{\cc_{\p^\perp}}\cdot e_{B} &= &\displaystyle\sum_{C\in cl(H)}e_{B}
(\widehat{C})\displaystyle\sum_{h\in C} {\cc}_{\p^\perp}\cdot h\\
{} & = & \displaystyle\sum_{C\in cl(H)}e_{B}(\widehat{C})
\displaystyle\sum_{h\in C}{\cc}_{(\p^\perp)^h},\\
{} & = & \displaystyle\sum_{C\in cl(H)}e_{B}(\widehat{C})
\displaystyle\sum_{h\in C}\sum_{\q\in(\p^{\perp})^h\cap I}{\cc}_{\q},

\end{array}
\end{equation*}
we have
\begin{equation*}
{\cc_{\p^\perp}}\cdot e_B =\sum_{\q\in I}\mathcal{S}_1(B,\p,\q){\cc}_{\q},
\end{equation*}
where
\begin{equation*}
\begin{array}{lll}
\mathcal{S}_1(B,\p,\q):=\displaystyle\sum_{C\in cl(H)}|
\mathcal{H}_{\p,\q}\cap C|e_B(\widehat{C}).
\end{array}
\end{equation*}
Similarly 
$\cc_{N(\p)}\cdot e_B=\sum_{\q\in I}\mathcal{S}_2(B,\p,\q)\cc_\q$
and 
$\cc_{\overline{N(\p)}}\cdot e_B=\sum_{\q \in I}\mathcal{S}_3(B,\p,\q)\cc_\q$,
where 
$$\mathcal{S}_2(B,\p,\q)=\sum_{C\in Cl(H)}|\mathcal{U}_{\p,N(\q)}\cap C|e_B(\widehat{C})$$
and
$$\mathcal{S}_3(B,\p,\q)=\sum_{C\in Cl(H)}|\mathcal{U}_{\p,\overline{N(\q)}}\cap C|e_B(\widehat{C}).$$

Assume first that $q\equiv 1\pmod 4$. 
If $\ell_{\p,\q}\in Pa_{\p}$, then $S_1(B_s,\p,\q)=0$ for each 
$s$ since $|\mathcal{H}_{\p,\q}\cap C|=0$ in $F$ for each 
$C\not=[0]$ by Lemma~\ref{m2}(i) and $e_{B_s}(\widehat{[0]})=0$ 
by Lemmas~\ref{expression} 2(c); and by Lemma~\ref{m2}(i) and 
Lemma~\ref{expression} 1(a), 1(c), 1(d), 3(a), 3(c), 3(d), 
we obtain
$$S_2(B_0,\p,\q)=e_{B_0}(\widehat{[0]})+e_{B_0}(\widehat{[\pi_k]})
+e_{B_0}(\widehat{D})=0+1+1=0$$
and
$$S_2(B_t^{'},\p,\q)=e_{B_t^{'}}(\widehat{[0]})+e_{B_t^{'}}
(\widehat{[\pi_k]})+ e_{B_t^{'}}(\widehat{D})=0+0+0=0.$$

If $\ell_{\p,\q}\in Se_{\p}$ and $\q \notin \p^\perp$, then 
by Lemma~\ref{m1}(ii) and Lemma~\ref{expression} 2(c) we obtain
$$S_1(B_s,\p,\q)=e_{B_s}(\widehat{[0]})+e_{B_s}(\widehat{[\theta_{i_1}]})
+e_{B_s}(\widehat{[\theta_{i_1}]})=0+0+0=0;$$
and by Lemma~\ref{expression} 1(c), 3(c), and Lemma~\ref{m2}(ii), 
$S_2(B_0,\p,\q)=e_{B_0}(\widehat{[0]})=0$ and
$S_2(B_t^{'},\p,\q)=e_{B_t^{'}}(\widehat{[0]})=0$.


If $\ell_{\p,\q}\in Se_{\p}$ and $\q\in \p^\perp$, 
then by Lemma~\ref{m1}(iii) and
Lemma~\ref{expression} 2(a) and 2(c) we
obtain
$S_1(B_s,\p,\q)=e_{B_s}(\widehat{[0]})+e_{B_s}(\widehat{D})=0+0=0$;
and from Lemma~\ref{m2}(ii) and 
Lemma~\ref{expression} 1(c) and 3(c), it follows that
$S_2(B_0,\p,\q)=e_{B_0}(\widehat{[0]})=0$
and
$S_2(B_t^{'},\p,\q)=e_{B_t^{'}}(\widehat{[0]})=0$.

Next we assume that $q\equiv 3\pmod 4$. 
If $\ell_{\p,\q}\in Pa_{\p}$ and $\q\notin\p^\perp$,
then by Lemma~\ref{m1}(v) and 
Lemma~\ref{expression} 5(c), we have 
$$S_1(B_r,\p,\q)=e_{B_r}(\widehat{[0]})+e_{B_r}
(\widehat{[\pi_{k_1}]})+e_{B_r}(\widehat{[\pi_{k_2}]})=0+0+0=0;$$
and by Lemma~\ref{m2}(iii) and Lemma~\ref{expression} 4(d) and 6(d), 
we obtain $S_3(B_0,\p,\q)=e_{B_0}(\widehat{[0]})=0$ and
$S_3(B_t^{'},\p,\q)=e_{B_t^{'}}(\widehat{[0]})=0$.

If $\q=\ell_{\p,\q}\cap\p^\perp$, then by Lemmas~\ref{m2}(iii)
and ~\ref{m1}(iii), and 4(d), 5(a), 5(c), 6(d) of Lemma~\ref{expression}, 
$S_3(B_0,\p,\q)= e_{B_0}(\widehat{[0]})= 0$,
$S_1(B_r,\p,\q) = e_{B_r}(\widehat{[0]})+
e_{B_r}(\widehat{D})= 0 + 0 = 0$, and 
$S_3(B_t^{'},\p,\q) = e_{B_t^{'}}(\widehat{[0]}) = 0.$

If $\ell_{\p,\q}\in Se_{\p}$, then by Lemmas~\ref{m2}(iv) 
and ~\ref{m1}(iv), and 4(a), 4(c), 4(d), 5(c), 6(a), 6(c), 
6(d) of Lemma~\ref{expression},
$$S_3(B_0,\p,\q) = e_{B_0}(\widehat{[0]})+e_{B_0}(\widehat{D})+
e_{B_0}(\widehat{[\theta_i]}) = 0 + 1 + 1 =0,$$
$S_1(B_r,\p,\q) = e_{B_r}(\widehat{[0]})=0$, and
$$S_3(B_t^{'},\p,\q) = e_{B_t^{'}}(\widehat{[0]}) + e_{B_t^{'}}
(\widehat{D})+e_{B_t^{'}}(\widehat{[\theta_i]})= 0 + 0 + 0 =0.$$
\QED

\noindent{{\bf Proof of Theorem~\ref{main}:}}  
Let $B$ be a $2$-block of defect $0$ of $H$. Then by 
Lemma~\ref{reduction}, we have
$$F^I\cdot e_B=\overline{\mathbf{S}^I\cdot f_B}.$$
Therefore, by Corollary~\ref{char11}, $F^I\cdot e_B=N$, 
where $N$ is the simple $FH$-module of dimension $q-1$ 
or $q+1$ lying in $B$ accordingly as $q\equiv 1\pmod 4$
or $q\equiv 3\pmod 4$.

Assume that $q\equiv 1\pmod 4$ and $q-1=m2^n$ with 
$2\nmid m$. Since
$$1=e_{B_0}+\displaystyle\sum_{s=1}^{(q-1)/4}e_{B_s}+
\displaystyle\sum_{t=1}^{(m-1)/2}e_{B_t^{'}},$$
$\Ker(\phi)\cdot e_{B_0}={\bf 0}$ and 
$\Ker(\phi)\cdot e_{B_t^{'}}={\bf 0}$,
then
$$\Ker(\phi)  = \displaystyle\bigoplus_{B\in Bl(H)}
\Ker(\phi)\cdot e_{B}= 
\displaystyle\bigoplus_{s=1}^{(q-1)/4}\Ker(\phi)\cdot e_{B_s} 
= \displaystyle\bigoplus_{s=1}^{(q-1)/4}N_s,$$
where $N_s$ is the simple module of dimension $q-1$ 
lying in $B_s$ for each $s$ by the discussion in the 
first paragraph. 

Now assume that $q\equiv 3\pmod 4$. Lemma~\ref{separate} 
yields 
$\Ker(\phi)=\langle\hat{\mathbf{J}}\rangle\oplus \Ima(\mu_2)$. 
Since $\Ima(\mu_2)\cdot e_{B_0}={\bf 0}$
and $\Ima(\mu_2)\cdot e_{B_t^{'}}={\bf 0}$,
applying the same argument as above, we have
$$\Ima(\mu_2)=\displaystyle\bigoplus_{r=1}^{(q-3)/4}M_r,$$
where each $M_r$ is a simple $FH$-module of dimension $q+1$.
Consequently, 
$$\Ker(\phi)=\langle\hat{\mathbf{J}}\rangle\oplus 
(\displaystyle\bigoplus_{r=1}^{(q-3)/4}M_r).$$
\QED

Now Conjecture~\ref{conj} follows as a corollary.

\begin{Corollary}
Let $\mathcal{L}$ be the $\Ff_2$-null space of 
$\A$. Then
\begin{equation*}
\dim_{\Ff_2}(\mathcal{L})=
\frac{(q-1)^2}{4}.
\end{equation*}
\end{Corollary}
{\Proof} By Theorem~\ref{main} and the fact that 
$\dim_{\Ff_2}(\mathcal{L})=\dim_{\Ff_2}(\Ker(\phi))$, 
when $q\equiv 1\pmod 4$, we have 
$$\dim_{\Ff_2}(\mathcal{L})=\displaystyle\sum_{i=1}^{(q-1)/4}(q-1),$$
and when $q\equiv 3\pmod 4$, we have
$$\dim_{\Ff_2}(\mathcal{L})=1+\displaystyle\sum_{i=1}^{(q-3)/4}(q+1),$$
both of which are equal to $\frac{(q-1)^2}{4}$.
\QED
\newpage
\section*{APPENDIX}
The character tables of $\PSL(2,q)$ were obtained by Jordan \cite{jordan}
and Schur \cite{schur} independently, from which we can
deduce the character tables of $H$ as follows.
Let $\epsilon\in \Cc$ be a primitive $(q-1)$-th root of
unity and $\delta\in \Cc$ a primitive $(q+1)$-th root of 
unity.

\begin{table}[htp]
\caption{Character table of $H$ when $q\equiv 1\pmod 4$}
\bigskip
\begin{tabular}{|c|c|c|c|c|c|c|c|}
\hline
 Number & 1 & 2 & $\frac{q-5}{4}$ & 1 & $\frac{q-1}{4}$ \\
\hline
 Size & 1 & $\frac{q^2-1}{2}$ & $q(q+1)$ & $\frac{q(q+1)}{2}$ & $q(q-1)$ \\
\hline
Representative &  $D$ & $F^{\pm}$ & $[\theta_i]$ & $[0]$ & $[\pi_k]$  \\
\hline
$\phi_{r} $ & $q+1$ & 1 & $\epsilon^{(2i)r}+\epsilon^{-(2i)r}$ & $2(-1)^{r}$& 0 \\

$\gamma$ & $q$ & 0 & 1& 1 & $-1$ \\

$1$ & 1 & 1 & 1 & 1 & 1 \\

$\chi_s$ & $q-1$ & $-1$ & 0 & 0 & $-\delta^{(2k)s}-\delta^{-(2k)s}$ \\

$\beta_1$ & $\frac{q+1}{2}$ & $\frac{1}{2}(1\pm\sqrt{q})$ & $\zeta(\theta_i)$ & $(-1)^{(q-1)/4}$ & 0 \\

$\beta_2$ & $\frac{q+1}{2}$ & $\frac{1}{2}(1\mp \sqrt{q})$ & $\zeta(\theta_i)$ & $(-1)^{(q-1)/4}$ & 0 \\

\hline
\end{tabular}\label{table:tab6.2}
\end{table}

Here $s=1,2,...,\frac{q-1}{4}$, $r=1,2,...,\frac{q-5}{4}$, $k=1,2,...,\frac{q-1}{4}$, $i=1,2,...,\frac{q-5}{4}$, and $\zeta(\theta_i) =1$ or $-1$. 

\begin{table}[htp]
\caption{Character table of $H$ when $q\equiv 3\pmod 4$}
\bigskip
\begin{tabular}{|c|c|c|c|c|c|c|c|}
\hline
 Number & 1 & 2 & $\frac{q-3}{4}$ & 1 & $\frac{q-3}{4}$ \\
\hline
 Size & 1 & $\frac{q^2-1}{2}$ & $q(q+1)$ & $\frac{q(q-1)}{2}$ & $q(q-1)$ \\
\hline
Representative &  $D$ & $F^{\pm}$ & $[\theta_i]$ & $[0]$ & $[\pi_k]$  \\
\hline
$\phi_r$ & $q+1$ & 1 & $\epsilon^{(2i)r}+\epsilon^{-(2i)r}$ & 0 & 0 \\

$\gamma$ & $q$ & 0 & 1 & $-1$ & $-1$ \\

$1$ & 1 & 1 & 1 & 1 & 1 \\

$\chi_s$ & $q-1$ & $-1$ & 0 & $-2(-1)^{s}$ &$-\delta^{(2k)s}-\delta^{-(2k)s}$ \\

$\eta_1$ & $\frac{q-1}{2}$ & $\frac{1}{2}(-1\pm\sqrt{-q})$ & $0$ & $(-1)^{(q+5)/4}$
&$-\zeta(\pi_k)$ \\

$\eta_2$ & $\frac{q-1}{2}$ & $\frac{1}{2}(-1\mp\sqrt{-q})$ & $0$ & $(-1)^{(q+5)/4}$
&$-\zeta(\pi_k)$ \\
\hline
\end{tabular}\label{table:tab6.4}
\end{table}

Here $s=1,2,...,\frac{q-3}{4}$, $r=1,2,...,\frac{q-3}{4}$, $k=1,2,...,\frac{q-3}{4}$, $i=1,2,...,\frac{q-3}{4}$, and $\zeta(\pi_k) =1$ or $-1$.


\end{document}